\documentclass[a4paper,12pt,leqno,twoside]{article}
\usepackage{color}
\usepackage{bm}
\usepackage{exscale}
\usepackage{amsmath}
\usepackage{amsfonts}
\usepackage{pagecolor,lipsum}
\usepackage{xcolor}
\usepackage{textcomp}
\usepackage{wasysym}
\usepackage{stmaryrd}
\usepackage{amscd}
\usepackage{graphicx}
\usepackage{amsxtra}
\usepackage{amssymb}
\usepackage{theorem}
 \usepackage{relsize}
\usepackage[final]{epsfig}
\usepackage{eqnarray}
\usepackage{verbatim}
\newtheorem{proposition}{Proposition}[subsection]
\newtheorem{definition}[proposition]{Definition}

\newtheorem{lemma}[proposition]{Lemma}
{\theorembodyfont{\rmfamily}\newtheorem{remark}[proposition]{Remark}}

\newtheorem{corollary}[proposition]{Corollary}

{\theorembodyfont{\rmfamily}}

\newfont{\abc}{cmtt10 scaled 1200}

\def\R{\mathbb{R}}

\def\Z{\mathbb{Z}}
\def\Q{\mathbb{Q}}
\def\U{\mathbb{U}}
\def\P{\mathbb{P}}

\def\P{\mathbb{P}}

\def\U{\mathbb{U}}

\def\X{\mathbb{X}}

\def\I{\mathbb{I}}
\def\ve{\varepsilon}

\def\ra{\rightarrow}
\def\cs{\symbol{35}}
\def\p{\partial}
\def\qed{\hfill $\Box$ \\}
\def\mm{\mbox}
\def\v{= \emptyset}
\def\n{\neq \emptyset}

\def\bp{\langle A \rangle}

\def\bp{\langle A \rangle}

\begin{document}
 \vspace*{-0.3cm}

\begin{center}\Large{\bf{The Higher Dimensional Positive Mass Theorem II}}\\
\smallskip
{\small{by}}\\
\smallskip
\large{\bf{Joachim Lohkamp}}
\end{center}

\vspace{0.5cm}
\noindent Mathematisches Institut, Universit\"at M\"unster, Einsteinstrasse 62, Germany\\
 {\small{\emph{e-mail:  j.lohkamp@uni-muenster.de}}}
\vspace{0.5cm}

 {\small{\center \tableofcontents}

{\contentsline {subsection}{\numberline {}References}{. .}}}

\medskip

\setcounter{section}{1}
\renewcommand{\thesubsection}{\thesection}
\subsection{Introduction} \label{introduction}
\medskip

In this paper, we prove the space-time positive mass conjecture in arbitrary dimensions. The proof generalizes the argument we have given to establish the Riemannian positive mass conjecture in [L1], which may also be described as the time-symmetric case of the general conjecture considered in the present paper.\\

The problem underlying this conjecture was already recorded in Pauli's 1921 book. In [P],Ch.61, the interpretation of the energy density of the gravitational field, contributing to the total mass, and its possibly negative sign, are explicitly described as debated issues, even before Einstein finalized his field equations [Ei] in 1916.  The positive mass conjecture then asserts that, nevertheless, the total mass for any non-vacuous isolated relativistic gravitational systems should be positive.\\

The interest in the higher dimensional form of general relativity evolved already in the 20ties, for instance, in attempts to incorporate all physical forces in one common framework. This, but also its weighty role in differential geometry, made the positive mass conjecture interesting in all dimensions.\\

\textbf{Known Cases} \, Until now, the positive mass conjecture could be established for spaces of dimension $\le 7$ and for spin manifolds.\\

In the late 70ties, Schoen, Yau [SY1]-[SY3] introduced an approach to treat this problem in $3$ dimensions using minimal and marginally outer trapped (hyper)surfaces. It was further developed by Eichmair, Huang, Lee and  Schoen [E1],[S],[EHLS] to cover dimensions $\le 7$.  The central problem with this strategy is that, in dimensions $>7$, these surfaces may be singular. By classical means, the impact of the singularities could not be controlled restricting the usability of this technique to low dimensions.\\

The only known alternative argument, which applies to a fairly broad class of spaces, is due to Witten [W],[PT]. It works for spin manifolds and it uses the Lichnerowicz formula for the associated Dirac operator to derive the conjecture. But the specifically spin geometric tools do not cover the still considerably wider non-spin case.\\

\textbf{General Case} \,  To derive the general space-time positive mass conjecture (Theorem 1), we employ skin structures on minimal and trapped surfaces to merge the hypersurface approach into a broader strategy.\\

 It starts with a reduction of the problem to the non-existence of isolated appearances of $S_{em}>0$. Here, $S_{em}$ is the energy-momentum scalar curvature,  we define as a counterpart to scalar curvature $S$, for initial data sets. Concretely, we show that $S_{em}>0$-islands cannot exist (Theorem 2). This is parallel to the reduction of the Riemannian conjecture to the non-existence of $S>0$-islands. The underlying deformation arguments may be of independent interest in the study of energy conditions.\\

This reduction rules out technicalities from the use of asymptotically flat geometry and allows us to transform the problem to a geometric one on compact manifolds. This way we can better focus on the core issue,  the presence of singular trapped surfaces once we leave the low dimensional special case.  We can now use skin structures and hyperbolic unfoldings to eliminate the singularities from surgery operations. As in the initial reduction step,  this can be organized to match the Riemannian case in [L1].\\

In the following sections of this chapter, we give an overview of the contents of the paper. In Ch.\ref{sh0} we introduce some basic concepts. Then, in Ch.\ref{sh}, we properly state our main results and  in  Ch.\ref{int2} we outline how to derive them.

\subsubsection{Total Mass and Constraint Equations} \label{sh0}
\bigskip

In this section we recall some basic concepts and explain the setup for the positive mass conjecture. We start with a space-time  $(X^{n+1}, g_X)$ satisfying Einstein's field equations in the form
\[G(g)_{\mu \nu}:= Ric(g)_{\mu \nu} - {\textstyle \frac{1}{2}}\cdot S(g) \cdot g_{\mu \nu} =  T_{\mu \nu},\]  where $Ric$ is the Ricci curvature, $S$ is the scalar curvature, $T$  the stress-energy tensor.\\

Now we \emph{assume} that there is a Cauchy hypersurface. This is a space-like hypersurface $(M^n,g) \subset (X^{n+1}, g_X)$ intersecting any inextensible, non-spacelike curve once.\\

This situation can be interpreted as the existence of some initial-data which determine the space-time from solving of Einstein field equations. General space-times need not to admit Cauchy hypersurfaces. But their existence is very plausible in the case of isolated systems, like galaxies, where we also expect that a Cauchy hypersurface will approximate the Euclidean one, near infinity.\\

Besides its topology and Riemannian metric, a Cauchy hypersurface carries also the information about its embedding, its second fundamental form $h$. Such a set of data is abstracted in the following concept, cf.Ch.\ref{ami0},A for some technical details.\\

\textbf{Definition 1} {\itshape \,  A triple $(M^n, g,h)$  is called an  \textbf{initial data set}, when $M$ is a smooth manifold, $g$ a $C^{2, \alpha}$-regular Riemannian metric and $h$ a $C^{1, \alpha}$-regular symmetric  $(0,2)$-tensor on  $M$, for some $\alpha>0$.}\\

Let us return to the motivating case where $(M^n, g,h)$ is a Cauchy hypersurface in a space-time  $(X^{n+1}, g_X)$ solving the field equations. Then we find compatibility relations between $T(g_X)$, $g$ and $h$, on $M^n$, the \textbf{Einstein constraint equations} for the local energy density $\mu:=G(g)_{\nu \nu}$, and  the local
momentum density, the 1-form $J(\cdot ):=G(g)(\p/\p \nu, \cdot )$, where, in this case, $\nu$ is the future oriented normal of $M^n$.\\

Employing the Gauss-Codazzi equations, we rewrite $\mu$ and $J$ on $(M^n, g,h)$ intrinsically in terms of $g$, $h$ and their derivatives. This gives the constraint equations, we write, for an abstract initial data set,  as definitions of $\mu$ and $J$.\\

\textbf{Definition 2} {\itshape \, For an initial data set $(M,g, h)$, we define
\begin{itemize}
  \item the \textbf{local energy density} $\mu:= \left(S(g)  - ( | h |_g^2 - (tr_g h )^2)\right)/2$
  \item the \textbf{local momentum density} $J := div_g h -   d(tr_g h)$
\end{itemize}}

The two quantities $\mu$ and $J$ return to their role as constraints when we locally develop a Cauchy hypersurface into
a solution of the field equations. Such a \emph{local} development exists due to work of Choquet-Bruhat [Cb].\\

Being back on $(X^{n+1}, g_X)$, it is physically plausible to expect the validity of the so-called dominant energy condition for $T$:  for any observer the local energy density appears non-negative and the local energy flow vector is non-spacelike, that is, not faster than light. Note that this does, however, not directly constrain the energy density of the gravitational field, also contributing to the total energy.\\

Again, this energy condition, admits an intrinsic description on $(M,g, h)$. \\

\textbf{Definition 3} {\itshape \, An initial data set $(M,g, h)$ satisfies the \textbf{dominant energy condition (DEC)} provided}
\[\mu \ge |J|_g.\]

For the more specific case of an isolated gravitational system, we consider a Cauchy hypersurface approximating a Euclidean geometry in the sense of an embedded manifold approaching a totally geodesic Euclidean hyperplane in  Minkowski space. Thus this concept includes both $g$ and $h$:\\

\textbf{Definition 4} {\itshape \, An initial data set $(M,g, h)$ is \textbf{asymptotically flat}, if it can be decomposed into a compact set $K\subset M$, the \textbf{core} of $M$, and an \textbf{end} $M\setminus K$, which admits a diffeomorphism $M\setminus K\cong \R^n\setminus B$, for some closed ball $B\subset \R^n$ such that
\[g-\delta \in   C^{2,p}_{-q}(M) \mm{ and } h \in  C^{1,\alpha}_{-q-1} (M)\]
where $\delta$ is a smooth symmetric $(0,2)$-tensor which equals $g_{Eucl}$ on $M\setminus K \cong\R^n\setminus B$ and
$q \in ((n-2)/2,n-2)$, $\alpha \in (0,1)$.}\\

To simplify the exposition, we usually ignore to write the diffeomorphic identifications and we think (and write) $\R^n\setminus B$ directly as a subset of $M$.
The weighted norm  for $ C^{k,\alpha}$-regular functions is defined as as the  $ C^{k,\alpha}$-norm on (a small neighborhood of) $K$ and on
$\R^n\setminus B$ as:
{\small \[ \| f\|_{C^{k,\alpha}_{-q}(\R^n\setminus B)} := \sum_{|I| \le k} \sup_x \left| |x|^{ | I| + q }(\partial^{I} f)(x) \right| +  \sum_{|I| = k} \left[ |x|^{\alpha + | I| + q }(\partial^{I} f)(x) \right]_\alpha\]}
Writing $f \in C^{k,\alpha}_{-q}$ means $\| f\|_{C^{k,\alpha}_{-q}(\R^n\setminus B)}<\infty$. In terms of Landau symbols this is written $f=O^{k, \alpha}(|x|^{-q})$. For  $f \in C^0_{-q}$ we write
 the more common $f = O(|x|^{-q})$.\\

(In more technical discussions, e.g. in [EHLS], one also considers weaker decay conditions in terms of weighted Sobolev norms and adds assumptions for $\mu$ and $J$.)\\

 Originally, Einstein [Ei], and others, [P],Ch.61, defined the total energy as volume integrals of energy contributions from the matter and the gravitational field over an asymptotically flat Cauchy hypersurfaces.   Later, these integrals were replaced by more intrinsically defined surface integrals over the sphere at infinity, using some Stokes formula argument. The following definition is due to Arnowitt, Deser and Misner[ADM].\\

\textbf{Definition 5} {\itshape \,  To an asymptotically flat initial data set $(M^n,g,h)$, we assign the following quantities, one oftentimes finds supplemented by the prefices \textbf{total} or \textbf{ADM}:
\begin{itemize}
\item the  \textbf{energy} $E (M,g)$, depending only on $g$,
\item the  \textbf{momentum} $P (M,g,h)=(P_1,...,P_n)$ of length $|P\,|= (\sum_i P_i^2)^{1/2}$,
\item the \textbf{energy-momentum}\, $\widehat{P} (M,g,h)=(E, P_1,...,P_n)$,
\item the  \textbf{mass}*\, $\bm{m}(M,g,h):=\sqrt{E (M,g)^2-|P (M,g,h)|^2}$,
\end{itemize}
 defined as follows:
{\small \begin{equation}\label{adm}
E (M,g) = \alpha_n\cdot \lim_{R \ra
\infty} \int_{\p B_R} \sum_{i,j} \left(\frac{\p g_{ij}}{\p x_i} - \frac{\p g_{ij}}{\p x_j} \right) \cdot \nu_j \: dV_{n-1},
\end{equation}
\begin{equation}\label{admp}
P_i (M,g,h) = 2 \cdot \alpha_n \cdot \lim_{R \ra
\infty} \int_{\p B_R} \sum_j \pi_{ij} \cdot \nu_j \: dV_{n-1}, \, i=1,...,n
\end{equation}}
where $\nu = (\nu_1\ldots \nu_n)$ is the outer normal vector to $\p B_R$ and  $\alpha_n= (2  (n-1) Vol(S^{n-1}))^{-1}$.}\\

*The mass can be thought of as the norm of $\widehat{P}$, with respect to the Lorentz metric, on the space-time development $X$ of $M$.
It is the validity of the positive mass conjecture, claiming $E \ge |P\,| \ge 0$,  which ensures that $\bm{m}(M,g,h)$ is a proper non-negative number.

\subsubsection{The Space-Time Positive Mass Theorem} \label{sh}
\bigskip

Now we can  state our main results, valid in any dimension $n \ge 3$.\\

\textbf{Theorem 1 (Space-Time Positive Mass Theorem)}{\itshape \,  Let $(M^n,g,h)$ be an asymptotically flat initial data set satisfying the dominant energy
condition. Then, we have \[E \ge |P\,| \ge 0.\]}
\emph{In Lorentzian geometric terminology, this says the total energy-momentum vector  $\bm{P}$  is a future directed time-like or null vector.}\\

This result settles other well-known versions of the positive mass conjecture, like the so-called space-time positive \emph{energy} conjecture, which only asserted $E \ge 0$. Also, for $h \equiv 0$, we reach the time-symmetric or Riemannian case, where $P=0$. Thus, the Riemannian positive \emph{mass} theorem just means $E \ge 0$.\\

Theorem 1 can be reduced to a basic obstruction result we formulate in Theorem 2 below. It applies to a variant $S_{em}(g,h)$ of the usual scalar curvature $S(g)$ adapted to the case of initial data sets.\\

\textbf{Definition 6}\, \emph{The \textbf{energy-momentum scalar curvature} $S_{em}(g,h)$ of an initial data set $(M^n,g,h)$ is the difference of the energy and momentum densities:}
 \begin{equation}\label{ems}
 S_{em}(g,h) := 2 \cdot(\mu - |J|_g)= S(g)  - ( | h |_g^2 - (tr_g h )^2) -  2 \cdot |div_g h  -  d(tr_g h)|_g
 \end{equation}

\smallskip

That is, the DEC can be expressed as $S_{em}(g,h) \ge 0$. And, in the time-symmetric or Riemannian case,  where $h \equiv 0$, we have $S_{em}(g,0) \equiv S(g) \equiv 2 \cdot \mu$. In this case, the DEC becomes $S\ge 0$. Now we claim there is \emph{no general local mechanism} to deform arbitrary initial data sets to others satisfying the strict form of the DEC. \\

\textbf{Theorem 2} \textbf{(Non-Existence of $\boldsymbol{S_{em} > 0}$-Islands) }{\itshape \,There exists \textbf{no} asymptotic flat initial data set $(M^n, g,h)$  such that:
 \begin{itemize}
   \item $S_{em}(g,h) >0$ on some non-empty open set $U \subset M^n$, with compact closure.
   \item $(M^n \setminus U, g,h) \equiv (\R^n \setminus  B_1(0), g_{Eucl},0)$.
 \end{itemize}}

The identification refers to the diffeomorphism from Def.4, which, in this case, is an isometry and this means that, on $M^n \setminus U$, $g$ is flat and $h \equiv 0$.\\

\subsubsection{Outline of Arguments and Techniques} \label{int2}
\bigskip

We start with the reduction of Theorem 1 to Theorem 2. It is based on deformation arguments employing geometric interpretations of $E$ and $P$. From this, Theorem 2 becomes the main result of this paper.
 This  matches the logic of the argument for the Riemannian positive mass theorem in [L1]. In that case, the non-existence of $S> 0$-islands was the corresponding main asssertion. \\

Assuming we had an $S_{em} > 0$-island, that is, a counterexample to Theorem 2, we can build particular compact  $S_{em} > 0$-manifolds. They are studied by means of marginally outer trapped surfaces and skin structures. We use the control this analysis gives to infer that such $S_{em} > 0$-manifolds cannot really exist. This contradiction establishes Theorem 2 and, therefore,  Theorem 1.\\

In the outline that follows, we describe this workflow and give some first ideas about the meaning and use of marginally outer trapped surfaces. From this we can also explain how to
approach the core problem in this business, the occurrence of singularities.\\

 \textbf{Reduction to $\boldsymbol{S_{em}>0}$-Islands} \, We assume there were a counterexample to the inequality $E \ge |P\,| \ge 0$ of Theorem 1, that is, a space with $E < |P\,|$. Indeed, we may always assume that $E < 0 \le  |P\,|$.\\

This is possible due to the boost argument of Christodoulou and  O'Murchadha [CO]. In simple terms, it implies we can deform the initial data set within some auxiliary space-time development, so that the new initial data set carries a Lorentz transformed energy-momentum vector $\widehat{P}=(E,P)$.  When we initially had  $0\le E < |P\,|, $ the transformation can be chosen, so that new vector has $E<0$.\\

Now, for $E < 0 \le  |P\,|$, we can apply deformations  trading negative energy for some positive $S_{em}$ and flatness towards
infinity. If $|P| >0$, the resulting space is an \emph{$S_{em} > 0$-peak}, this is an initial data set $(\P^n, g_\P,h_\P)$ with $S_{em}>0$ such that:
 \[(\P^n \setminus U, g_\P,h_\P) \equiv (\R^n \setminus  B_R(0),w^{4/n-2} \cdot g_{Eucl}, w^{2/n-2} \cdot h^\triangle),\mm{ for some } R>0,\]
 where  $w=1-\alpha \cdot |x|^{-(n-2)} -  \beta \cdot  |x|^{-(n-1)}$, for some $\alpha,\beta >0$, and $h^\triangle$ denotes the following symmetric form derived from the Euclidean Green's function $1/|x|^{n-2}$:
 {\small \[h^\triangle_{nn}=\frac{2 \cdot \p (|x|^{-(n-2)})}{\p x_n},\, h^\triangle_{nj}= h^\triangle_{jn}=\frac{\p (|x|^{-(n-2)})}{\p x_j}\mm{ and  } h^\triangle_{ij}=0 \mm{ for } i,j < n.\]}
When $P=0$, we can even choose $w=1$ and $h \equiv 0$ on $M^n \setminus U$. That is, in this case, we directly get an  $S_{em} > 0$-island.\\

 To eliminate $h^\triangle$, occurring in the $|P| >0$-case, we exploit its particular shape to superpose two copies of an $S_{em} > 0$-peak, but with momentum vectors pointing in opposite directions. The resulting space has $E<0$, $P=0$  with $S_{em}>0$. Thus we can now use this \emph{twin peaks space} to again build an $S_{em} > 0$-island. \\

\textbf{Trapped Surfaces} \, To disprove the existence of such an $S_{em} > 0$-island, that is, a counterexample to Theorem 2, we employ \emph{marginally outer trapped surfaces}, abbreviated MOTS. To explain this concept, we consider a possibly singular submanifold $H^{n-1}$, of an initial data set $(M^n,g,h)$, with second fundamental form $h$ relative $M^n$.\\

$H^{n-1}$ is called a MOTS if its mean curvature, relative $M^n$, annihilates the mean curvature contributions of $M^n$: $\theta^+_H:=tr_H h_{M|H} +  tr_H h_H=0$. Here,
 $h_{M|H}$ is the restriction of $h$ on the ambient space $M^n$ to vectors tangent to $H^{n-1}$, while
$tr_H h_H$ denotes the trace of the second fundamental form $h_H$ of $H^{n-1} \subset M^n$.\\

To better understand this condition, let $L^{n+1}$ be a space-time development of  $M^n \subset L^{n+1}$. Then, being a MOTS means that the future volume expansion $\theta^+_H$, in outward null-direction, vanishes. The typical context is that of black holes. Then, MOTS are not properly trapped  surfaces, with $\theta^\pm < 0$, eventually collapsing in a space-time singularity. Instead, they are limiting or marginal, since the outer expansion vanishes $\theta^+=0$, making MOTS models for apparent horizons of black holes.\\

MOTS are not known to be representable as minimizers of an elliptic variation problems. Nevertheless, there is a notion of \emph{stable} MOTS  sharing many analytic properties with stable, and even with area minimizing, hypersurfaces. Formally, they are almost minimizers and this also means they generally admit singularities in dimensions $\ge 7$. The concept and the existence of stable MOTS became available through work of Andersson, Mars, Simon, Metzger and Eichmair in [AMS1],[AM],[E2] and [EM].\\

\textbf{MOTS and $\boldsymbol{S_{em}}>0$} \, Stable MOTS in a $S_{em} > 0$-space, can be conformally deformed into $S>0$-geometries. This is the counterpart of the well-known $S>0$-heredity principle saying a stable minimal hypersurface in a $S > 0$-space can be conformally deformed into a $S>0$-geometry.\\

This idea to use stable MOTS evolved from a paper by Galloway and Schoen [GS] treating horizons of black holes.  It was implemented to derive Theorem 1 in dimensions $\le 7$, in the work of Eichmair, Huang, Lee and  Schoen [E1] and [EHLS]. They follow a classical strategy to inductively select, analyze and deform asymptotically flat stable MOTS in an asymptotically flat ambience, from [S].\\

Our setup is different and it simplifies the work with MOTS, in particular, in the hard case of singular MOTS. We transplant the hypothetical $S_{em} > 0$-island onto some large and flat torus. After some arbitrarily $C^3$-small perturbation, this gives us a $S_{em} > 0$-manifold of the form $\Q^n=T^n \cs N^n$, for some compact manifold $N^n$.\\

More specifically, this means $S>0$ and $h \equiv 0$ outside the core. Also whenever we choose two parallel tori $T_i^{n-1} \subset T^n$, $i=1,2$, relative to the flat initial metric, not intersecting the core, we can easily ensure that, now with respect to such an $S> 0$-metric,  they are \emph{strictly mean convex} when viewed as the two boundary components of the domain which contains the core.\\

 To prove Theorem 2, we show that such a geometric configuration cannot exist.\\

\textbf{Singular MOTS} \, The geometry on $\Q^n$ allows us to find a stable  MOTS $Y^{n-1} \subset \Q^n=T^n \cs N^n$ homologous to $T^{n-1} \subset T^n \cs N^n$, and lying between the two tori $T_i^{n-1}$.\\

 The regularity theory for these MOTS shows that $Y^{n-1}$ contains a nearly flat and large torus component: $Y^{n-1}$ can again be written $Y^{n-1}= T^{n-1} \cs N^{n-1}$, for some compact but, in dimensions $>7$, usually \emph{singular} space $N^{n-1}$.\\

With the occurrence of these singularities, we reach the principal problem which classically obstructed the use of minimal and marginally outer trapped (hyper)surface techniques, in dimensions $>7$.
For a discussion of some of the resulting and quite peculiar issues we refer to [L4],Ch.1.1.C and Ch.1.3.A.\\

At this stage, we appeal to skin structural techniques from [L2]-[L4]. We may apply them to \emph{compact} singular MOTS equipped with their stability operators and conformal Laplacians. There are several reasons why the genuine skin structural arguments do \emph{not} apply to \emph{asymptotically flat} singular MOTS. This is one of the essential benefits from starting the argument with a compactification via $S_{em} > 0$-islands.\\

\textbf{Skin Structures and MOTS} \, Skin structures are used on different levels. We first observe that the compact stable MOTS  $Y^{n-1}$, very much like area minimizing hypersurfaces, is skin uniform. This says the singular set $\Sigma_Y$, viewed as a boundary of $Y \setminus \Sigma_Y$, has a boundary regularity (precisely) sufficient  to gain a  detailed control over the asymptotic analysis of many elliptic operators towards $\Sigma_Y$ on the original space $Y$.\\

This control is invested in a further and differently natured application of skin structures, an \emph{$S_{em}>0$-to-$S>0$-heredity principle with surgeries}.\\

The outcome is a that for any neighborhood $V$ of $\Sigma_{N^{n-1}}$, we can conformally deform $Y^{n-1}$ into some space $Z^{n-1}$  so that there is another, smoothly bounded, neighborhood $U \subset V$ of $\Sigma_{Z^{n-1}}=\Sigma_{Y^{n-1}}$ so that, with respect to the new metric:
\begin{itemize}
\item $\Q^{n-1}=Z^{n-1} \setminus U$ has $scal >0$ and $\p \Q^{n-1}=\p U$ has positive mean curvature.
\end{itemize}

Also, this deformation reproduces the second geometric main property of $T^n \cs N^n$:
 \begin{itemize}
\item $\Q^{n-1}$ also contains a nearly flat and large torus component.
\end{itemize}

But this sort of space cannot exist. This follows from a largely parallel argument of $S>0$-heredity employing (singular) minimal hypersurfaces. Indeed, at this stage, we simply switch to the, from now onwards, identical steps in proof of the Riemannian positive mass conjecture in [L1]. The outcome is an inductive reduction scheme finally reaching a  \emph{non-existing}  $S>0$-surface of genus $\ge 1$.\\

This contradiction completes the proof of Theorem 2 and, thus, of Theorem 1.\\

\setcounter{section}{2}
\renewcommand{\thesubsection}{\thesection}
\subsection{Reduction to $S_{em}>0$-Islands} \label{ssosm}
\bigskip

Given a counterexample to Theorem 1, that is, a space  $(M^n, g, h)$, satisfying the DEC, but with $E < |P\,|$, we first get $S_{em} > 0$-peaks and then $S_{em} > 0$-islands, that is, counterexamples to Theorem 2.\\

In Ch.\ref{ami0} we introduce a collection of tools to deform initial data sets while keeping the DEC. Some of them are new and may also be useful in other contexts.
Then, in the main section Ch.\ref{ami2}, we apply them to get $S_{em} > 0$-islands we use in the next chapter to build particular geometries amenable to a study by stable MOTS.

\subsubsection{Basic Deformations} \label{ami0}
\bigskip

In this section we collect some auxiliary deformations for initial data sets before we start with the actual reduction to $S_{em} > 0$-islands and compact models.\\

\textbf{A. Conformal Deformations} \, To modify $S_{em}$, we  conformally deform the metric $g$ of an initial data set $(M^n,g,h)$, to $v^{4/n-2} \cdot g$, for some $C^2$-function $v>0$.\\

Now, we recollect the transformation rules for the scalar curvature $S$,  the covariant derivative $\nabla^v$ and for $h$ under this conformal change,cf.[Be],Ch.1J:
\begin{equation}\label{sc}
 \textstyle S(v^{4/n-2} \cdot g) \cdot v^{\frac{n+2}{n-2}} =  -\gamma_n \cdot \Delta v + S(g) \cdot v,\, \mm{ with } \gamma_n = \frac{4 (n-1)}{n-2},
\end{equation}
\begin{equation}\label{cvd}
 \textstyle  \nabla^v_X Y = \nabla_X Y+ \frac{2}{n-2} \cdot  \left((X \log v) \cdot Y+ (Y \log v) \cdot X - \langle X,Y \rangle \cdot \nabla \log v \right).
\end{equation}
\begin{equation}\label{h01}
|v^{2/n-2} \cdot  h |_{v^{4/n-2} \cdot g}^2=v^{-4/n-2} \cdot | h |_g^2  \, \mm{ and }\,  (tr_{v^{4/n-2} \cdot g} v^{2/n-2} \cdot  h )^2=v^{-4/n-2} \cdot (tr_g h )^2
\end{equation}

\begin{remark} \,\label{me} For the two identities of (\ref{h01}), we recall that the form $h$ is assumed to be an arbitrary $(0,2)$-tensor on $M^n$. However, what is not made explicit in the literature, is the expected transformation behavior of $h$, under changes of $g$. It is defined to match the equality $h \cdot \nu = A$, where $A$ is the mean curvature vector and $\nu$ the unit normal vector to $M^n \subset X^{n+1}$ in some space-time development. That is, conformal changes of $g$, thought as restrictions from $X$, also affect the length of $\nu$.\\

To compensate this additional scaling effect, one considers the conformal deformed initial data sets of the form $(M^n,v^{4/n-2} \cdot g, v^{2/n-2} \cdot h)$. The term $v^{2/n-2} \cdot h$ occurs because, abstractly, $(M^n,v^{4/n-2} \cdot g)$ does not arise from an embedding. However, if $M^n \subset L^{n+1}$ was a Cauchy hypersurface and we conformally deform $L$, with gradient parallel to $M$, then the new second fundamental form $h_{new}$ for $M$ is $h_{new}=v^{2/n-2} \cdot h$. \qed
\end{remark}

Also, for a local orthonormal frame field $e_1,...e_n$, we get for the $1$-form $J$
\begin{equation}\label{j}
 J_i=J(e_i)=\sum_{j=1} (\nabla_{e_j} h)(e_i,e_j) - D_{e_i} h(e_j,e_j)
\end{equation}
where $(\nabla_X h)(Y,Z):=D_X h(Y,Z) - h(\nabla_X Y,Z) - h(Y, \nabla_XZ) $ for vector fields $X,Y,Z$. \\

Thus, to understand the changes of $J$, under conformal transformations $(M^n, g, h)$ to $(M^n,v^{4/n-2} \cdot g, v^{2/n-2} \cdot h)$, it suffices to consider the two emerging types of terms:
\begin{equation}\label{jees}
D_{v^{-2/n-2} \cdot e_i} (v^{2/n-2} \cdot h)(v^{-2/n-2} \cdot e_j,v^{-2/n-2} \cdot e_k) =
\end{equation}
\[v^{-2/n-2} \cdot \left(v^{-2/n-2} \cdot \p  h_{jk}/ \p x_i + \p v^{-2/n-2}/\p x_i  \cdot h_{jk}\right),\]
\begin{equation}\label{jas}
v^{2/n-2} \cdot h(\nabla^v_{v^{-2/n-2} \cdot e_i}{v^{-2/n-2} \cdot e_j},{v^{-2/n-2} \cdot e_k})=
\end{equation}
\[ v^{-2/n-2} \cdot  \left(\p v^{-2/n-2}/\p x_i   \cdot h_{jk}+ v^{-2/n-2} \cdot h(\nabla^v_{e_i}{ e_j},e_k)\right).\]
Also, for the  transformation of the norm $|J|= (\sum_i^n J_i^2)^{1/2}$ of $J$, we get:
\begin{equation}\label{trj}
|J(a^{4/n-2} \cdot g,a^{2/n-2} \cdot h)|_{a^{4/n-2} \cdot g} = a^{-4/n-2} \cdot |J( g,h)|_g,
\end{equation}
for any constant $a >0$.\\

\textbf{B. Local $\boldsymbol{S_{em}}>0$-Deformations} \, We want to locally redistribute or scatter $S_{em}>0$ using conformal deformations  $(M^n,v^{4/n-2} \cdot g, v^{2/n-2} \cdot h)$. To this end, we use some specific  cut-off functions $\phi=f(d,s)$, for  $s,d>0$:
\begin{equation}\label{fds}
f (d,s) (r) := s \cdot exp (-d/r) \mbox{ on } \mathbb{R}^{\ge 0}, f (d,s) \equiv 0 \mbox{ on } \mathbb{R}^{\le 0}.
\end{equation}
One readily checks that for any $0<\ve_2 < \ve_1<1$, there is a  $d=d(\ve_1,\ve_2)\gg 1$  with
\begin{equation}\label{fe}
0 < f(d,s) < \ve_1 \cdot f'(d,s) < \ve_2 \cdot f''(d,s) \mm{ on } (0,10).
\end{equation}
Also, we choose some fixed cut-off function for some $\chi \in C^\infty(\R,\R^{\ge 0})$, with $\chi \equiv 1$ on $\R^{\le 0}$ and
 $\chi \equiv 0$ on $\R^{\ge 1}$ and define $\psi(d,s) (r) := \chi(r-8) \cdot  f (d,s) (r).$\\

We use these functions as follows: for some smooth diffeomorphism $A_V:(-1,10) \times S^{n-1} \ra V$ onto an annulus $V \subset M$, we define the function $\Psi(d,s)$ on $V$: \[\Psi(d,s)(x):=\psi(d,s)(\pi_1 \circ A^{-1}_V(x))\]
where $\pi_1$ is the projection on the first coordinate. Then we have

\begin{lemma}\label{jesq}  \, For any $\kappa \ge 1$, there is a $d_0(A_V,g,\kappa)\gg 1$ so that for  $d \ge d_0, s>0$:
\begin{equation}\label{fee}
|\Delta \Psi(d,s)| \ge \kappa \cdot \max\{|\nabla \Psi(d,s)|, |\Psi(d,s)|\}, \mm{ on } A_V((0,8) \times S^{n-1}).
\end{equation}
\end{lemma}

\textbf{Proof} \,  We write  the Laplacian of $\Psi$ in terms of coordinates $x_1,x_2..x_n$ we get from  $(-1,10) \times S^{n-1}$ via $A_V$:  $\Delta \, \Psi = \frac{1}{\sqrt{det \, g}} \cdot \sum_{i,j} \frac{\p}{\p x_i} (\sqrt{det \, g} \cdot g^{ij} \cdot \frac{\p \, \Psi}{\p x_j})$. (Formally, we cover $ S^{n-1}$ by finitely many local coordinate patches and get for each of them coordinates $x_2,..,x_n$, whereas $x_1$ is the canonical coordinate for $(-1,10)$.)\\

In this expression, the only non-vanishing second  derivative is $\frac{\p^2 \, \Psi}{\p x_1^2}$. On compact subsets of $V$, we observe that its coefficient $c_2(A_V,g,x)=g^{11}(x)$, in the positive definite matrix $(g^{ij})$, is positively lower bounded and the coefficient $c_1(A_V,g,x)$ for the only non-trivial first derivative
$\frac{\p \, \Psi}{\p x_1}$ is bounded.\\

 Thus we infer from  (\ref{fe}) that, for large $d$, and any compact subset $K \subset V$, $x \in K$:
 \[\textstyle |\Delta \Psi(d,s)(x)| = \Big|c_2(A_V,g,x) \cdot \frac{\p^2 \, \Psi}{\p x_1^2}  + c_1(A_V,g,x) \cdot \frac{\p \, \Psi}{\p x_1}\Big|\ge a_1(A_V,K,g) \cdot \frac{\p^2 \, \Psi}{\p x_1^2},\] for some  $a_1(A_V,K,g)>0$. Similarly, we get an $a_2(A_V,K,g)>0$, so that  $|\nabla \Psi(d,s)| \le a_2(A_V,K,g) \cdot |\frac{\p \, \Psi}{\p x_1}|$.  Another application  of (\ref{fe})  gives the estimate (\ref{fee}).\qed

Now we  write, using  (\ref{sc}):
 \begin{equation}\label{sse}
S_{em}(v^{4/n-2} \cdot g,  v^{2/n-2} \cdot h) = (-\gamma_n \cdot \Delta v + S(g) \cdot v) \cdot v^{-(n+2)/(n-2)}...
 \end{equation}
  \[...- (v^{-\frac{4}{n-2}} \cdot | h |_g^2  - v^{-\frac{4}{n-2}} \cdot (tr_g h )^2)  -  |J(v^{\frac{4}{n-2}} \cdot g, v^{2/n-2} \cdot h)|_{v^{\frac{4}{n-2}} \cdot g}.\]

As a consequence of (\ref{cvd})-(\ref{trj}), saying that all perturbation terms are upper estimated by bounded multiples of $\max\{|\nabla \Psi(d,s)|, |\Psi(d,s)|\}$,  and (\ref{fee}) we observe:

\begin{corollary}\label{jesc}  \,  If $S_{em}(g,h) \ge 0$ on $V$ and $S_{em}(g,h) > 0$ on $A_V((7,10) \times S^{n-1})$,
then we have, for large $d$ and sufficiently small $s>0$:
\begin{equation}\label{def} S_{em}((1-\Psi(d,s))^{4/n-2} \cdot g,(1-\Psi(d,s))^{2/n-2} \cdot h)>0 \mm{ on } A_V((0,10) \times S^{n-1}).
\end{equation}
\end{corollary}
In typical applications of this result we have $S_{em}(g,h) > 0$ on an open subset $U \subset M$ and $S_{em}(g,h) \ge 0$ outside.  Then \ref{jesc} allows us to scatter $S_{em} > 0$ to any connected
superset of $U$. Somewhat more generally we have:

\begin{corollary}\label{jesc2}  \,  Assume that $S_{em}(g,h) > 0$ on $U \subset M$ and $S_{em}(g,h) \ge 0$ outside. Also let $W \subset M$ be open with compact closure $K_W \cap   \overline{U} \v$, and $(M,g_1,h_1)$ be an initial data set, with $g_1\equiv g$ and $h_1\equiv h$, outside $W$.\\

Then, for any open and connected set $U^* \subset M$ with $U \cup K \subset U^*$, there is an $\ve(U,U^*,K_W,g,h) >0$, so that if $|g_1-g|_{C^2} \le \ve$ and  $|h_1-h|_{C^2} \le \ve$, there is conformal deformation  $v^{4/(n-2)} \cdot g$, with
\[v\equiv 1\mm{ outside } U^* \, \mm{ and } \, S_{em}(v^{4/(n-2)} \cdot g_1,v^{2/(n-2)} \cdot h_1) > 0 \mm{ on }U^*.\]
\end{corollary}
\textbf{Proof} \,  We first consider the case $\ve =0$. Since $U^* \subset M$ is open and connected, we can find a (possibly finite) sequence of open sets $X_i$,$i \ge 1$ with $\bigcup_ i X_i=U^*$ so that there are diffeomorphisms of Euclidean balls $\iota_i: B_3(0) \ra X_i$, with
\begin{itemize}
  \item  $\iota_1(B_2(0)) \subset U$ and $\iota_{i+1}(B_2(0))  \subset U \cup \bigcup_{1 \le k \le i} X_k$,
  \item for any compact $K \subset U^*$, there is an $i_K$, so that for $i \ge i_K$:   $X_i \cap K \v$.
\end{itemize}

 Then we can inductively apply \ref{jesc} to scatter $S_{em}>0$ from $U$ to $U \cup \bigcup_{1 \le k \le i} X_k$, for any $i\ge 1$, without accumulations. For $i \ra \infty$, we get the assertion for $U^*$.\\

For $\ve >0$ small enough, we notice that  perturbations of $g$ and $h$ supported on $U^*$, of $C^2$-norm $\le \ve$, do not change the validity of the open condition $S_{em}>0$. \qed

\textbf{C. Harmonic Asymptotics} \, Many global properties of asymptotically flat initial data sets depend only on the decay of some key terms we derive from in $g$ and $h$. To exploit this observation, we want to give  $g$ and $h$ a more transparent look. To this end, we approximate the original data by more elementary ones, keeping the main properties unchanged. \\

 Here we employ the density result of Eichmair, Huang, Lee and  Schoen, [EHLS], Th.18.. This perturbation result  says that, while keeping $E <  |P\,|$ and the DEC, we can assume that $(M, g, h)$ has the following type of asymptotics.

\begin{definition}\label{hay}  \,  An asymptotically flat initial data set $(M^n,g,h)$ has \textbf{harmonic asymptotics} if there exist a $C^{2,\alpha}$-function $u$ and a $C^{2,\alpha}$-vector
field $Y=(Y_1,...Y_n)$ on $M$ with
\begin{equation}\label{11}
 u(x) = 1+ a \cdot |x|^{-(n-2)}+ O^{2, \alpha}(|x|^{-(n-1)}) \mm{ and } Y_i(x)=b_i \cdot |x|^{-(n-2)}+O^{2, \alpha}(|x|^{-(n-1)}),
\end{equation}
\begin{equation}\label{22}
 g=u^{\frac{4}{n-2}}\cdot g_{Eucl} \mm{ and } h_{ij}=u^{\frac{2}{n-2}}\cdot (Y_{i,j}+Y_{j,i}), \mm{ outside some compact } K \subset M
\end{equation}
 for $i,j= 1, ..., n$ and where $a$, $b_1, ...b_n$ are constants.
\end{definition}

The most interesting terms in the development of $u$ and $Y_i$ are Green's (and this means harmonic) functions and one has:
\begin{equation}\label{ep}
\textstyle a=\frac{1}{2} \cdot E \, \mm{ and }\,  b_i = - \frac{n-1}{n-2} \cdot P_i  \mm{ for } i= 1, ..., n,
\end{equation}
for the energy $E$ and momentum $P$ of $(M^n,g,h)$, cf.[EHLS], Lemma 5. A nice aspect of the density theorem [EHLS], Th.18 is its flexibility towards the DEC.
One may alternatively choose $(M^n,g,h)$ satisfying the strict DEC, $\mu > |J|$, or so that $\mu = |J|=0$ near infinity, simulating the vacuum case. We will later use both options in succession to readjust the shape of our initial data set.\\

\textbf{D. Geometric Asymptotics} \,  In our reduction of Theorem 1 to 2 we transform the relation between total energy and momentum into a geometric constraint. For this we deform harmonic asymptotics near infinity into geometrically more appealing ones.\\

 The deformations are arranged in layers and raise the decay towards infinity. Up to higher order perturbations we consider the transitions from $r^\beta$ to $r^\alpha$ to $\alpha<\beta<0$, where we allow $\alpha =-\infty$ and set $r^{-\infty}=0$ on $\R^{\ge 2}$. We recall that for any  $H \in C^2(\R^{>0},\R)$, $\Delta H(|x|) = H'' (r) + \frac{n-1}{r} \cdot H' (r)$, $r=|x|$, on $\R^n \setminus \{0\}$. For $|x|^\alpha$ we have
\begin{equation}\label{de}
\Delta |x|^\alpha= \big(\alpha \cdot (\alpha -1) +  (n-1) \cdot \alpha \big)\cdot |x|^{\alpha-2}>0, \mm{ for } \alpha<-(n-2),
\end{equation}
and for the Green's function $G_n(x)=|x|^{-(n-2)}$ we get $\Delta |x|^{-(n-2)}=0$.

\begin{lemma} \label{flay}\emph{\textbf{(Subharmonic Deformations towards  Infinity)}}  \\

\textbf{A. Decay Switches} \,  For $\kappa \ge 1$ and  $m(\kappa,n)\ge 1$ large enough, there is large $d(\kappa,n)>0$ and an $s_0(d) >0$ and an $a(n,d,s)> 0$, so that for $s \in (0,s_0(d)]$, there are smooth functions $F[Z,\kappa]_{n,m,d,s}$, for $Z=A,B$, on $\R^+$ with
\begin{enumerate}
\item $F[Z,\kappa]>0$, $(r^{-(n-2)})' \le F[Z,\kappa]' \le 0$, $F[Z,\kappa]'' \ge 0$, $\Delta F[Z](|x|)  \ge 0$, $r=|x|$,
\item $F[Z,\kappa](r)=r^{-(n-2)}+f(d,s \cdot (m+2)^{-(n-2)})(r-(m-1))$ on $(0,m]$,
\item  $\Delta F[Z,\kappa](|x|)  \ge a/m^{n-1}$ on $B_{m+3}(0) \setminus \overline{B_m(0)}$.
\end{enumerate}
On $\R^{\ge m+3}$ we have for $Z=A,B$
{\small
\begin{equation}\label{fab}
F[A](r):= F[A,1](r)=(m+2)^{-(n-2)} +f\big(d,s/(m+2)^{n-2} \big)((m+4)-r),
\end{equation}
\[F[B,\kappa](r)= (m+2)^{-(n-2)} + \kappa \cdot \big(r^{-(n-1)} - (m+3)^{-(n-1)}\big) +f\big(d,s/(m+2)^{n-2} \big)((m+4)-r).\]}
\textbf{B. Perturbed Green's Functions} \, For $G_n+ f,$ for some $f=O^{2, \alpha}(|x|^{-(n-1)})$, and any sufficiently large $m$, there is a smooth function $F[C]=F[C][m,f]>0$ with:
\begin{enumerate}
\item $F[C]=G_n+ f$ on $\overline{B_{m+1}(0)}$ and $F[C] \equiv (m+2)^{-(n-2)}$ on  $\R^n \setminus  B_{m+3}(0)$,
\item $|\nabla F[C]| \le 2 \cdot |\nabla G_n|$ on $B_{m+5}(0) \setminus \overline{B_{m}(0)}$,
\item $\Delta F[C] \ge \Delta f$ on $B_{m+2}(0) \setminus \overline{B_{m+1}(0)}$ and $\Delta F[C] > 0$ on $B_{m+3}(0) \setminus \overline{B_{m+2}(0)}$.
\end{enumerate}
\end{lemma}

The proof is a combination of elementary cut-off and smoothing constructions very similar to those in [L5],Lemma 6.2. We leave the adaptation of the details to the reader.\qed

In our applications, we employ combinations of these functions with coefficients we choose depending on the total energy $E$ and momentum $P$ of our initial data set:
\begin{equation}\label{abc}
 v_{\Theta,\kappa}:=1 -\theta_0 \cdot G_n(x)- \theta_A \cdot  F[A](|x|) - \theta_B \cdot F[B,\kappa](|x|) -  \theta_C \cdot   F[C](x),
\end{equation}
for some $\Theta=(\theta_0, \theta_A,\theta_B,\theta_C)$, with $\theta_A=0$, $\theta_0,\theta_B,\theta_C >0$, when $P \neq 0$ and $ \theta_A,\theta_C>0$, $\theta_0=\theta_B=0$, when $P = 0$. For conformal deformations of an asymptotically flat space, we have from (\ref{sc})  and $S(g_{Eucl})=0$:
{\small \begin{equation}\label{ac}
\textstyle S(v_{\Theta,\kappa}^{4/n-2} \cdot g) \cdot v_{\Theta,\kappa}^{n+2/n-2}= \frac{4 (n-1)}{n-2} \cdot \big(\theta_A \cdot \Delta F[A](|x|) +\theta_B \cdot \Delta F[B,\kappa](|x|) + \theta_C \cdot  \Delta F[C](x)\big).
\end{equation}}

\textbf{E. Boost to }$E < 0 \le  |P\,|$ \, Now we want to show that we may assume that $E < 0 \le  |P\,|$. This essentially is a classical (folklore) argument explained
e.g. in [EHLS], p.119, 3rd Remark and its references.\\

From subsection C, we may assume that  $(M^n,g,h)$ has  harmonic asymptotics and satisfies the DEC with $\mu = |J|=0$ near infinity. Then the boost techniques of Christodoulou and  O'Murchadha [CO] show that there is space-time development of  $(M^n, g, h)$, solving the field equations, for a period of time linearly increasingly with the distance from the core, cf. [CO], Th.6.1. \\

Within this so-called boost domain, we can deform the initial data set to another asymptotically flat one still  satisfying the DEC. The deformation asymptotically becomes a Lorentz transform of the end.  The point is that this also Lorentz transforms its energy-momentum vector $\widehat{P}=(E, P)$.\\

In our case,  [CO], Th.6.1, shows that the boost domain asymptotically has slope $1$, that is, eventually exhausts the complement of a light cone. Therefore, we can be sure that when we initially have $0 \le E < |P\,|$, that is, if the energy-momentum vector is space-like, we can use a boost to Lorentz transform  $\widehat{P}$ into another vector with $E<0$. Thus we can henceforth assume that $E < 0$.\\

From this boost argument,  applying the density theorem a second time and some scaling and rotation, we can henceforth assume that our initial data set  $(M^n,g,h)$ has a more specific shape.

\begin{definition}\label{res}  \, An asymptotically flat initial data set $(M^n,g,h)$ has  a \textbf{$\diamond$-reduced shape}*, provided it has the following properties:
\begin{itemize}
\item  an asymptotically flat end with  \textbf{harmonic asymptotics}.
\item \textbf{strictly negative energy}, $E =-4< 0 \le   |P\,|$,
 \item \textbf{strict DEC}, $\mu > |J|$, in geometric terms $S_{em}>0$,
\end{itemize}
Concretely, outside a compact $K \subset M$ we have
\begin{equation}\label{ex}
g= u^{4/n-2} \cdot g_{Eucl}, \mm{ for } u(x)  = 1-2\cdot |x|^{-(n-2)}+ f_u\mm{ and }f_u=O^{2, \alpha}(|x|^{-(n-1)})
\end{equation}
Depending on the value of the total momentum $P$, we have the cases
\begin{itemize}
\item $\boldsymbol{P =0}$ \, $b_1=...b_n=0$ and this means
\begin{equation}\label{h10}
h_{ij}= O^{1, \alpha}(|x|^{-n})
\end{equation}
\item $\boldsymbol{P \neq 0}$\, $P=(0,...0, - \frac{n-2}{n-1})$ and, thus, $b_1=...b_{n-1}=0$ and $b_n=1$. That is
\begin{equation}\label{h0}
 h_{ij}= u^{\frac{2}{n-2}}\cdot h^\triangle_{ij} + O^{1, \alpha}(|x|^{-n})=  O^{1, \alpha}(|x|^{-(n-1)})
\end{equation}
where $h^\triangle$ is the primitive form we get from the Green's function $1/|x|^{n-2}$
 {\small \begin{equation}\label{h1}
h^\triangle_{nn}=\frac{2 \cdot \p (|x|^{-(n-2)})}{\p x_n},\, h^\triangle_{nj}= h^\triangle_{jn}=\frac{\p (|x|^{-(n-2)})}{\p x_j}\mm{ and  } h^\triangle_{ij}=0 \mm{ for } i,j < n.\end{equation}}
\end{itemize}
\end{definition}
(*The prefix $\diamond$ is meant to symbolize the distinctive inequalities $E<0, S_{em}>0$.)\\

\subsubsection{$E  < |P|$ versus $S_{em}>0$} \label{ami2}
\bigskip

We show that for a given asymptotically flat initial data set we can invest negative total energy to transform the given space into another one with  some changed end structure amenable to further geometric arguments. These modifications employ relations between the total energy and $S_{em}$, depending on the value of the total momentum.

\begin{definition}\emph{\textbf{($\boldsymbol{S_{em}>0}$-Islands and  Peaks)}}\label{pks} \, We consider asymptotic flat initial data sets  $(\I^n, g_\I,h_\I)$, $(\P^n, g_\P,h_\P)$ and $(M^n,g,h)$.
Then we use the following terminology.\\

$\bullet$ \, $(\I^n, g_\I,h_\I)$ is an $\boldsymbol{S_{em}>0}$\textbf{-island}, provided
\begin{enumerate}
\item $S_{em}(g_\I,h_\I) >0$ on some non-empty open set $U \subset \I^n$, with compact closure.
\item $(\I^n \setminus U, g_\I,h_\I) \equiv (\R^n \setminus  B_1(0), g_{Eucl},0)$.
\end{enumerate}

$\bullet$ \, $(\P^n, g_\P,h_\P)$  is an  $\boldsymbol{S_{em}>0}$\textbf{-peak} if it has $S_{em}(g_\P,h_\P)>0$  on $\P^n$ and
\[(\P^n \setminus U, g_\P,h_\P) \equiv (\R^n \setminus  B_R(0),w^{4/n-2} \cdot g_{Eucl}, w^{2/n-2} \cdot h^\triangle),\mm{ for some } R>0,\]
where  $w=1-\alpha \cdot |x|^{-(n-2)} -  \beta \cdot  |x|^{-(n-1)}$, for some $\alpha,\beta >0$.\\

$\bullet$ \, $\boldsymbol{C}_m:=M \setminus \left(\R^n \setminus B_m(0)\right)$ is the  \textbf{core} of radius $m \ge 10.$
\end{definition}

\begin{proposition}\emph{\textbf{(Reduction to $S_{em}>0$-Islands and Peaks)}} \label{redo}  We deform an initial data set $(M^n,g,h)$ of $\diamond$-reduced shape
into initial data sets $(M^n, g[m,\kappa],h[m,\kappa])$ with $g[m,\kappa]=g$ and $h[m,\kappa]=h$ on $\boldsymbol{C}_{m-1}$, and so that on $\R^n \setminus B_{m-1}(0)$ for $v_\Theta$ of (\ref{abc})
\begin{itemize}
\item $\boldsymbol{P=0} \quad  g[m,1]= v_{\Theta,1}^{4/n-2} \cdot g_{Eucl},  \,  h[m,1]=v_{\Theta,1}^{2/n-2} \cdot \omega_m \cdot h$, 
\[\mm{ for }\theta_A=\theta_C =1/2,\theta_0=\theta_B=0 \mm{ and } \kappa  = 1,\]
\item $\boldsymbol{P\neq 0} \quad   g[m,\kappa]= v_{\Theta,\kappa}^{4/n-2} \cdot g_{Eucl},  \, h[m,\kappa]= v_{\Theta,\kappa}^{2/n-2} \cdot \big(h^\triangle +  \omega_m \cdot (h-h^\triangle)\big)$, \[\mm{ for }\theta_A=0,\theta_0=\theta_B=1/4, \theta_C =1/2 \mm{ and } \kappa \ge 1,\]
\end{itemize}
and where $\omega_m(x):=\omega(x-(m+1))$, $x \in \R$, for some cut-off function $\omega \in C^\infty(\R,\R^{\ge 0})$, with $\omega \equiv 1$ on $\R^{\le 0}$ and
 $\omega \equiv 0$ on $\R^{\ge 1}$. Then, we have:
\begin{itemize}
\item $\boldsymbol{P=0} \quad   (M^n, g[m,1],h[m,1])$ is an $S_{em}>0$-island, for  large $m$,
\item $\boldsymbol{P \neq 0} \quad  (M^n, g[m,\kappa],h[m,\kappa])$ is an $S_{em}>0$-peak, for large $\kappa$ and $m$, where $m$ also depends on the chosen $\kappa$.
\end{itemize}
\end{proposition}

\textbf{Proof} \, We separate $(M^n, g[m,\kappa],h[m,\kappa])$ into 4 differently treated regions
\[M^n=[\boldsymbol{C}_{m-1}] \cup [\boldsymbol{C}_m\setminus  \boldsymbol{C}_{m-1}  \cup \boldsymbol{C}_{m+4}\setminus  \boldsymbol{C}_{m+3}] \cup [\boldsymbol{C}_{m+3}\setminus  \boldsymbol{C}_m] \cup [M^n \setminus \boldsymbol{C}_{m+4}].\]

In the $S_{em}>0$-peak case, we encounter when $P \neq 0$, we also need the parameter $\kappa$, and this happens precisely in the case of dimension $n=3$. On the first 3 regions, we get $S_{em}>0$, for any $\kappa \ge 1$, when we successively choose $m$ large enough. Only in the last step, on $[M^n \setminus \boldsymbol{C}_{m+4}]$, we also need to choose $\kappa \gg 1$ to ensure that $S_{em}>0$. \\

$[\boldsymbol{C}_{m-1}]$\, $g[m,\kappa]=g$ and $h[m,\kappa]=h$ and hence $S_{em}(g[m,\kappa],h[m,\kappa])>0$.\\

$[\boldsymbol{C}_m\setminus  \boldsymbol{C}_{m-1}  \cup \boldsymbol{C}_{m+4}\setminus  \boldsymbol{C}_{m+3}]$\,  In the case $P=0$, we add the terms
$- f(d,|E|/4 \cdot s)(r-(m-1))$ respectively $- f(d,|E|/4 \cdot s)((m+4)-r)$ of \ref{flay} to the original conformal factor, while we have $h[m,\kappa]=h$.  For large $d$ and small $s$, we infer,  from \ref{jesc}   that $S_{em}(g[m,\kappa],h[m,\kappa])>0 $.\\

For $P \neq 0$, the same argument applies on $\boldsymbol{C}_m\setminus  \boldsymbol{C}_{m-1}$, whereas on $\boldsymbol{C}_{m+4}\setminus  \boldsymbol{C}_{m+3}$ we extend the argument from $[M^n \setminus \boldsymbol{C}_{m+4}]$ we discuss below. \\

$[\boldsymbol{C}_{m+3}\setminus  \boldsymbol{C}_m]$\, We use the decay estimate  $\Delta F(Z,\kappa)(|x|)  \ge a/m^{n-1}$ from \ref{flay}.A(ii) and  (\ref{abc}) to positively lower estimate $S(g[m,\kappa])$. We show that all other $S_{em}(g[m,\kappa],h[m,\kappa])$-terms can be upper bounded by $O(|x|^{-n})$-functions, that is, by functions $f$ with
\begin{equation}\label{ewq}
|f| \le c_f \cdot  |x|^{-n},\mm{ for some } c_f >0, x \in \R^n \setminus  B.
\end{equation}
The  point is that the constant $c_f$, and the function $f$, can be chosen \emph{independently} of $m$.  From this we infer that, for  $m \gg 1$, $S(g[m,\kappa])$  exceeds these stronger decaying terms in  $S_{em}(g[m,\kappa],h[m,\kappa])$, on $\boldsymbol{C}_{m+3}\setminus  \boldsymbol{C}_m$.\\

In what follows, $x_1,..,x_n$ denote the Euclidean coordinates and $e_1,..e_n$ the associated orthonormal and parallel frames on $(\R^n \setminus  B, g_{Eucl})$. That is,
we   have $\langle e_i,e_j\rangle=\delta_{ij}$ and $\nabla_{e_i} e_j=0$, for any pair $i,j$. With respect to $g[m,\kappa]$, the $v_{\Theta,\kappa}^{-2/n-2} \cdot e_i$, $i=1,..n$, will be again an orthonormal, but generally non-parallel frame field. so that from (\ref{cvd}):
\begin{equation}\label{nagg}
\textstyle \nabla^{v_{\Theta,\kappa}}_{e_i} e_k =  \frac{2}{n-2} \cdot  \left((\p  v_{\Theta,\kappa}/\p x_i) \cdot e_k+ (\p v_{\Theta,\kappa}/\p x_k) \cdot e_i- \delta_{ij} \cdot \nabla  v_{\Theta,\kappa}\right)/v_{\Theta,\kappa}.
\end{equation}

\textbf{Estimates for  $v_{\Theta,\kappa}$}: For $a_{n,m}:=1-3/4 \cdot (m+2)^{-(n-2)}+ \kappa/4 \cdot (m+3)^{-(n-1)}=lim_{x \ra \infty}v_{\Theta,\kappa}(x)$. For $\alpha \in \R$, some $c_{\alpha,n} >0$, we get from a Taylor expansion of $z^\alpha$ in  $a_{n,m}$,
\begin{equation}\label{ba1}
|v_{\Theta,\kappa}^\alpha -a_{n,m}^\alpha|\le  c_{\alpha,n}\cdot |v_{\Theta,\kappa}-a_{n,m}| \le f_0 \mm{ and } |\nabla v_{\Theta,\kappa}| \le 2 \cdot |\nabla G_n| \le f_1
\end{equation}
for $ f_0= O(|x|^{-(n-2)})$, $ f_1= O(|x|^{-(n-1)}),$ both independent of $\kappa$ and $m$, once we choose $m$, depending on $\kappa$, large enough. From (\ref{nagg}) and (\ref{ba1})
\begin{equation}\label{nnn}
|D_{v_{\Theta,\kappa}^{-2/n-2} \cdot e_i} v_{\Theta,\kappa}^{-2/n-2}|, |\nabla^{v_{\Theta,\kappa}}_{v_{\Theta,\kappa}^{-2/n-2} \cdot e_i}{v_{\Theta,\kappa}^{-2/n-2} \cdot e_k}| \le  10 \cdot   |\nabla G_n|\le f_2,
\end{equation}
for $f_2 = O(|x|^{-(n-1)})$, in the same sense, independent of $\kappa$ and $m$. \\

\textbf{Estimates for $\mu$}: from (\ref{h10})-(\ref{h1}) and  (\ref{ba1}) we get
\begin{equation}\label{h}
|h[m,\kappa]|_{g[m,\kappa]}^2, (tr_{g[m,\kappa]} h[m,\kappa])^2  \le f_3,
\end{equation}
for some $f_3=  O(|x|^{-2 \cdot (n-1)}) $, as before, independent of $\kappa$ and $m$. \\

\textbf{Estimates for $J$}: To estimate the decay of $|J(g[m,\kappa],h[m,\kappa])|_{g[m,\kappa]}$, we consider the terms in (\ref{jees}) and (\ref{jas}). To get the shape assumed there, one introduces the factors $v_{\Theta,\kappa}^{2/n-2}$ for $h^\triangle$ and $(v_{\Theta,\kappa}/u)^{2/n-2}$ for $u^{2/n-2} \cdot h$. We note from (\ref{ba1}) that
 \begin{equation}\label{dev}
 (v_{\Theta,\kappa}^{2/n-2}-a_{n,m}^{2/n-2})  \cdot h^\triangle_{ij},\, \big((v_{\Theta,\kappa}/u)^{2/n-2}-a_{n,m}^{2/n-2}\big)   \cdot h_{ij} = O^{1, \alpha}(|x|^{-n}).
\end{equation}
 Thus, for the decay estimates we may absorb this deviation in the decay term and omit writing it explicitly. Then, we first observe from  (\ref{h10})-(\ref{h1})  and (\ref{nnn}) that
\begin{equation}\label{also}
|h[m,\kappa](\nabla^{v_{\Theta,\kappa}}_{e_i}{ e_j},e_k)| \le |h[m,,\kappa]| \cdot f_2 \le  f^*_0, \mm{ for } f^*_0 = O(|x|^{-2 \cdot (n-1)})
\end{equation}
For the other terms, we distinguish between
the cases $P  = 0$ or $P \neq 0$:\\

$\bullet \, \boldsymbol{P=0}$ \, We have, readily from (\ref{h10}):
{\small \[|\p \omega_m/\p x_k  \cdot h_{ij}|  \le f^*_4, \mm{ for } f^*_4 =  O(|x|^{-n}) \mm{ and } |\omega_m  \cdot \p h_{ij}/\p x_k|  \le f^*_5, \mm{ for } f^*_5 =   O(|x|^{-(n+1)})\]}

$\bullet \, \boldsymbol{P \neq 0}$ \, We get from (\ref{h0}) and  (\ref{h1}):
{\small \[|\p h^\triangle_{ij}/\p x_k| \le f^*_1, \mm{ for } f^*_1= O(|x|^{-n}) \mm{ and } |\omega_m \cdot \p (h_{ij} - h_{ij}^\triangle)/\p x_k| \le f^*_2,\mm{ for } f^*_2=  O(|x|^{-(n+1)})\]
\[\mm{ and } |\p \omega_m/\p x_k  \cdot (h_{ij} - h_{ij}^\triangle)| \le f^*_3,\mm{ for } f^*_3 =  O(|x|^{-n})\]}
This is the main difference between the cases $P = 0$ and $P \neq 0$. For $P \neq 0$ we do not get the $O(|x|^{-n})$ estimate for $h$ but only for $h-h^\triangle$, since $h^\triangle=O(|x|^{-(n-1)})$.\\

Thus we have in both cases, $P  = 0$ and $P \neq 0$,  from (\ref{ba1})
\begin{equation}\label{sum1}
\textstyle |\p  h[m,\kappa]_{jk}/ \p x_i| \le f^*_6  \mm{ and }  \p v^{\alpha}/\p x_i  \cdot h[m,\kappa]_{jk}\le f^*_7, \mm{ for } \alpha=-\frac{2}{n-2}, -\frac{4}{n-2},
\end{equation}
for  $f^*_6 =  O(|x|^{-n})$,  $f^*_7 = O(|x|^{-2 \cdot (n-1)})$. Using again  (\ref{ba1}), we also see that the terms from (\ref{also}) and (\ref{sum1}), multiplied by $v_{\Theta,\kappa}^{\beta}$, $\beta=-\frac{2}{n-2}, -\frac{4}{n-2}$, all satisfy these bounds. Summing up,   this shows that
 \begin{equation}\label{jjj}
 |J(g[m,\kappa],h[m,\kappa])|_{g[m,\kappa]} \le f^*_8,\mm{ for }f^*_8 =  O(|x|^{-n}),
 \end{equation}
where all $f^*_i$,$i=1,..,8$, are independent of $m$, where again,  we choose $m$, depending on $\kappa$, large enough\\

\textbf{Estimates for $S_{em}$}: From our assumption that $E<0$, and $S(g_{Eucl})=0$, we get for any given $\kappa\ge 1$ and sufficiently large $m$:
\[S_{em}(g[m,\kappa],h[m,\kappa])= -\gamma_n \cdot \Delta v_{\Theta,\kappa} \cdot v_{\Theta,\kappa}^{-\frac{4}{n-2}}  - \big[| h[m,\kappa] |_{g[m,\kappa]}^2  - (tr_{g[m,\kappa]} h[m,\kappa] )^2\big]-...\]
\[...- |J(g[m,\kappa],h[m,\kappa])|_{g[m,\kappa]} \ge  \frac{b_n}{m^{n-1}} - \left[2 \cdot \frac{c_{f_3}}{m^{2 \cdot (n-1)}} + \frac{c_{f^*_7}}{m^n} \right]>0, \mm{ on }  \boldsymbol{C}_{m+3}\setminus  \boldsymbol{C}_m,\]
for  some  $b_n>0$.\\

$[M^n \setminus \boldsymbol{C}_{m+4}]$\,  For  $P=0$ we have $g[m,\kappa]=g_{Eucl}$ and $h[m,\kappa]=0$. The non-trivial case is
$P \neq 0$. In this case, we claim that, for $n>3$ and $\kappa=1$ respectively, for $n=3$ and sufficiently  large $\kappa \gg 1$, there is a large $R_1>0$ so that
\begin{equation}\label{as}
S_{em}(v_{\Theta,\kappa}^{4/n-2} \cdot g_{Eucl}, v_{\Theta,\kappa}^{2/n-2} \cdot h^\triangle) >0 \mm{ on } \R^n \setminus B_{R_1}(0).
\end{equation}
To  check this, we note that since $S(g_{Eucl})=0$ and $\Delta |x|^{-(n-2)}=0$, the transformation law for $S$ under conformal deformations reads $S(v_{\Theta,\kappa}^{4/n-2} \cdot g_{Eucl}) \cdot v_{\Theta,\kappa}^{\frac{n+2}{n-2}} =  -\gamma_n \cdot \Delta v_{\Theta,\kappa}$ with
 {\small \begin{equation}\label{vv}
\Delta v_{\Theta,\kappa}=\frac{\p^2 v_{\Theta,\kappa}}{\p r^2} +\frac{n-1}{r} \cdot \frac{\p v_{\Theta,\kappa}}{\p r}=- \kappa \cdot \big(n \cdot (n-1) -  (n-1)^2 \big)\cdot |x|^{-(n+1)}/2<0.
 \end{equation}}
This shows, for any $\kappa \ge 1$, there is a large $R_0>0$
 \begin{equation}\label{ses}
 S(v_{\Theta,\kappa}^{4/n-2} \cdot g_{Eucl}) \ge  \kappa \cdot |x|^{-(n+1)} \mm{ on } \R^n \setminus B_{R_0}(0).
 \end{equation}
As in (\ref{h}) we get, using again  Taylor expansions of $z^\alpha$ in  $1$: \begin{equation}\label{qc}
\left||v_{\Theta,\kappa}^{2/n-2} \cdot h^\triangle|_{v_{\Theta,\kappa}^{4/n-2} \cdot g_{Eucl}}^2 - (tr_{v_{\Theta,\kappa}^{4/n-2} \cdot g_{Eucl}} v_{\Theta,\kappa}^{2/n-2} \cdot  h^\triangle)^2\right|  \le f[1],
\end{equation}
for some $f[1]=  O(|x|^{-2 \cdot (n-1)}) $. Also, computing $J$ in terms of the vector field $Y$ used in (\ref{22}) to define harmonic asymptotics with $h_{ij}= Y_{i,j}+Y_{j,i}$
\[\textstyle J_i = (div_{g_{Eucl}} h -   d(tr_{g_{Eucl}} h))(e_i)=\sum_j\Big(\p \big(\p Y_i/\p x_j + \p Y_j/\p x_i - \delta_{ij} \cdot \sum_k \p Y_k/\p x_k )\big)/\p x_j\Big)\]
From this we have, for $h^\triangle=O^{1, \alpha}(|x|^{-(n-1)})$, $Y=(0,...0, 1/|x|^{n-2})$, a $z^\alpha$-Taylor expansion, in (\ref{jees}) and (\ref{jas}), and the product rule for differentiations:
{\small \begin{equation}\label{jn}
 J(g_{Eucl},h^\triangle)_i=\Delta Y_i = 0, \mm{ since }  Y_i \equiv 0, \mm{ for } i < n , \, J(g_{Eucl},h^\triangle)_n=\Delta (1/|x|^{n-2}) = 0,
\end{equation}
\begin{equation}\label{q2}
|J(v_{\Theta,\kappa}^{4/n-2} \cdot g_{Eucl}, v_{\Theta,\kappa}^{2/n-2} \cdot h^\triangle)|_{v_{\Theta,\kappa}^{4/n-2} \cdot g_{Eucl}} \le
\end{equation}
\[O(|x|^{-(n-2)})  \cdot  O(|x|^{-n}) + O(|x|^{-(n-1)})  \cdot  O(|x|^{-(n-1)}) = :  f[2],\]}
so that $f[2]=  O(|x|^{2 \cdot (n-1)})$. The upper bounds (\ref{qc}) and (\ref{q2}) in terms of $f[i]$ apply only on $M^n \setminus \boldsymbol{C}_{m+4}$, but they do not depend on $m$. Thus, the exponents for the decay of the $S_{em}$-term and the other terms in $S_{em}$ satisfy the following relations
\begin{equation}\label{na}
n+1  < 2 \cdot (n-1) \mm{ for } n >3 \mm{ and } n+1=2 \cdot (n-1) \mm{ for } n =3.
\end{equation}
Therefore, for $ n >3$, we infer from (\ref{ses}) that, for $\kappa =1$ and some large $R_1 \ge R_0>0$, (\ref{as}) holds and this means  $S_{em}>0$ on $M^n \setminus \boldsymbol{C}_{m+4}$, for $m \ge R_1$.\\

For $ n =3$, we use \ref{flay} for some large $\kappa \gg 1$ and observe that the choice of $\kappa$ influences the radii $R_1,R_0$, and thus the size of $m$. But it does not affect the estimates (\ref{qc}) and (\ref{q2}) outside $\boldsymbol{C}_{m+4}$, for a sufficiently large $m$ since the Green's function term and its derivative eventually overcompensate all contributions from $\kappa \cdot |x|^{-(n-1)}$. Thus we also get  $S_{em}>0$ on $M^n \setminus \boldsymbol{C}_{m+4}$, for $ n =3$, when we choose $\kappa \gg 1$ and hence $m$ large enough.\qed

\begin{lemma}\emph{\textbf{(Twin Peaks and Islands)}}\label{pks2}   If there exists an $S_{em}>0$-peak, then we also get an $S_{em}>0$-island.
\end{lemma}
\medskip
\textbf{Proof} \quad  We first define a twin peaks space $\X^n$ from a superposition of two  $S_{em}>0$-peaks. Then we show that, placing the cores of the peaks far enough from each other,  we can ensure that $S_{em}>0$ and get an initial data set of $\diamond$-reduced shape with $P=0$.\\

\textbf{Twin Peaks} \, We choose some $\kappa$ and $m\ge 10$ large enough so that \ref{redo}  for $P \neq 0$ gives an $S_{em}>0$-peak.  Outside the core $\boldsymbol{C}_{m+5}$ we have
\begin{itemize}
 \item $g[m,\kappa]=v_{\Theta,\kappa}^{4/n-2} \cdot g_{Eucl}$ and $h[m,\kappa]=v_{\Theta,\kappa}^{2/n-2} \cdot  h^\triangle$,
\end{itemize}
for $v_{\Theta,\kappa}(x)=\psi_{m,\kappa}(r):=a_{n,m}- 1/2 \cdot  r^{-(n-2)} - \kappa/2 \cdot  r^{-(n-1)}$, $r=|x|$.\\

For $e_1:=(1,0,..,0) \in \R^n $, $D\ge 10^2 \cdot m+10^2$, the translation $T_{D \cdot e_1}(x)=(x+D \cdot e_1)$ we assemble an $S_{em}>0$-twin peaks space $\X^n$ from attaching two suitably oriented copies $\boldsymbol{C}^\pm_{m+5}$ of $\boldsymbol{C}_{m+5}$ disjointly  to $\R^n$. Topological, we have
\[\X^n=\big(\R^n \setminus  B_{m+5}(D \cdot e_1) \cup B_{m+5}(-D \cdot e_1)\big) \cup \boldsymbol{C}^+_{m+5} \cup \boldsymbol{C}^-_{m+5}.\]

To define an initial data set $(\X^n,g_D,h_D)$, we start with the following superposition of deformations on $\R^n \setminus  B_{m+5}(D \cdot e_1) \cup B_{m+5}(-D \cdot e_1)$, where a cut-off function $\varphi$ neutralizes the influence of the deformation, centered in $\pm D \cdot e_1$, in $B_{m+8}(\mp D \cdot e_1)$.
\begin{itemize}
\item $\Phi_{m,\kappa, D}(x):=\varphi_m(x-D \cdot e_1) \cdot \psi_{m,\kappa} (|x+D \cdot e_1|)+\varphi_m(x+D \cdot e_1) \cdot \psi_{m,\kappa} (|x-D \cdot e_1|)$,
\item $g_D=\Phi_{m,\kappa,D}(x)^{4/n-2} \cdot g_{Eucl}$,\, $h_D= \Phi_{m,\kappa,D}(x)^{2/n-2} \cdot  (T^*_{-D \cdot e_1}(h[m,\kappa]) -  T^*_{D \cdot e_1}(h[m,\kappa]))$,
\end{itemize}
with $\varphi_m(x)=\overline{\varphi}(|x|-m)$, for some smooth function $\overline{\varphi} \in C^\infty(\R,[0,1])$ with $\overline{\varphi}=0$ on $\R^{\le 8}$ and $\overline{\varphi}=1$ on $\R^{\ge 9}$.\\

For $\boldsymbol{C}^+_{m+5}$ we glue one copy of the core $\boldsymbol{C}_{m+5}$ along $\p B_{m+5}(D \cdot e_1)$, using the restriction of the identity map of $\R^n$.\\

For $\boldsymbol{C}^-_{m+5}$ we glue another copy of $\boldsymbol{C}_{m+5}$ along $\p B_{m+5}(-D \cdot e_1)$. But now this is done after a $180^\circ$-rotation of $\boldsymbol{C}_{m+5}$ in the $(x_1,x_n)$-plane.\\

Under a $180^\circ$-rotation in the $(x_1,x_n)$-plane, the vector field $Y= (0,..0,Y_n) = |x|^{-(n-2)} \cdot \p/\p x_n$, underlying the definition of $h^\triangle$, transforms into $-Y$. Hence, in both cases,  the metric and the form on $\boldsymbol{C}_{m+5}$ smoothly extend along the gluing boundary and define a smooth initial data set $(\X^n_D,g_D,h_D)$.\\

\textbf{$\boldsymbol{S_{em}>0}$ near Infinity} \,We define the larger cores  $\boldsymbol{C}^\pm_{m+5+l} \subset \X^n$, for $0 \le l  \le m$,
\[\boldsymbol{C}^\pm_{m+5+l} := \boldsymbol{C}^\pm_{m+5} \cup B_{m+5+l}(\pm D \cdot e_1) \setminus B_{m+5}(\pm D \cdot e_1).\]
We show that for $m$ large enough and any $D\ge 10^2 \cdot m+10^2$
\begin{equation}\label{as}
S_{em}(g_D,h_D)>0 \mm{ on } \R^n \setminus  B_{m+10}(D \cdot e_1) \cup B_{m+10}(-D \cdot e_1).
\end{equation}
 As in the argument for \ref{redo}, case $[M^n \setminus \boldsymbol{C}_{m+4}]$, we infer that for large $\kappa$ and $m$:
 \begin{equation}\label{aa}
 S(g_D) \ge \kappa \cdot ( |x-D \cdot e_1|^{-(n+1)}+ |x+D \cdot e_1|^{-(n+1)})
 \end{equation}
{\small \begin{equation}\label{aaa}\left||h_D|_{g_D}^2 - (tr_{g_D} h_D)^2\right|  \le \end{equation}
\[4 \cdot \big(|T^*_{-D \cdot e_1}(h^\triangle)|_{g_{Eucl}}^2+|T^*_{D \cdot e_1}(h^\triangle)|_{g_{Eucl}}^2 + (tr_{g_{Eucl}} T^*_{-D \cdot e_1}(h^\triangle))^2+ (tr_{g_{Eucl}} T^*_{D \cdot e_1}(h^\triangle))^2\big)\]
\[ \le C_1 \cdot (|x-D \cdot e_1|^{-2 \cdot (n-1)}+ |x+D \cdot e_1|^{-2 \cdot (n-1)}),\]}
for some $C_1 >0$ independent of $\kappa$ and $m$, when we choose $m$ large enough. Also, since $J$ is linear in $h$, we have from $J(g_{Eucl},h^\triangle)=0$
\[J\left(g_{Eucl},  T^*_{-D \cdot e_1}(h^\triangle) -  T^*_{D \cdot e_1}(h^\triangle)\right)=0.\]
This shows, as in (\ref{q2}), from $\psi_{m,\kappa}-a_{n,m} =O^{2, \alpha}(|x|^{-(n-2)})$ and $h^\triangle=O^{1, \alpha}(|x|^{-(n-1)})$  that
\begin{equation}\label{aaaa}
|J(g_D,h_D)|_{g_D} \le C_2 \cdot (|x-D \cdot e_1|^{-2 \cdot (n-1)}+ |x+D \cdot e_1|^{-2 \cdot (n-1)}),
\end{equation}
for some $C_2 >0$, independent of $\kappa$ and $m$, again for  $m$ large enough. From (\ref{na})  we have, for dimension $n >3$,  combining (\ref{aa}) - (\ref{aaaa}):\, $S_{em}(g_D,h_D)>0$ on $\R^n \setminus  B_{m+10}(D \cdot e_1) \cup B_{m+10}(-D \cdot e_1)$, for large enough $m$.  For $n=3$ we argue as after (\ref{na}) and infer that $S_{em}(g_D,h_D)>0$, for sufficiently large $\kappa$ and  $m$.\\

\textbf{Separating the Cores}  \, Now we turn to $(\boldsymbol{C}^\pm_{m+10},g_D,h_D)$. Extending the last step we get $S_{em}(g_D,h_D)>0$ on $\boldsymbol{C}^\pm_{m+10} \setminus \boldsymbol{C}^\pm_{m+9}$, for sufficiently large $m$.\\

For $\boldsymbol{C}^\pm_{m+9}$,  we choose a basepoint $p_0 \in \boldsymbol{C}_{m+5}$, before we glue it onto $\R^n$ to define $\X^n$, and the associated point $p_0(D) \in \boldsymbol{C}^+_{m+5} \subset \X^n_D$. Then, we observe that, for $D \ra \infty$, the pointed initial data sets $(\X^n_D,g_D,h_D;p_0(D))$  compactly $C^3$-converge to  $(\P, g_\P,h_\P; p_0)$ which has $S_{em}>0$ everywhere. Thus, for $D$ large enough we
 have $S_{em}(g_D,h_D)>0$ on $\boldsymbol{C}^+_{m+10}$. The same argument applies to  $\boldsymbol{C}^-_{m+10}$.\\

 \textbf{Energy and Momentum}  \, For such a large $D>0$, $(\X^n,g_D,h_D)$ is asymptotically flat with harmonic asymptotics and it has $S_{em}>0$. Moreover, we now observe that, for large enough $m$, $|a_{n,m}-1| \le 1/10$ and thus $E(\X^n,g_D,h_D)<-1$. Due to the twin peaks definition of $(\X^n,g_D,h_D)$  we get
 \begin{equation}\label{x}
(h_D)_{ij}= O^{1, \alpha}(|x|^{-n}) \mm{ and, in particular, } P(\X^n,g_D,h_D)=0
\end{equation}
 This follows from the observation that, for $t \in \R$ and  any given $a \in \R$:
{\small \begin{equation}\label{ki}
 \left|\frac{1}{t^{n-1}}-\frac{1}{(t+a)^{n-1}}\right| = \left|\frac{\sum_{i=1}^{n-1}\binom{{n-1}}{i} t^{{n-1}-i}\cdot a^i}{(t+a)^{n-1} \cdot t^{n-1}}\right| =O(t^{-n}).
\end{equation}}

Summarizing, we assembled an initial data set of $\diamond$-reduced shape with $P=0$. Now we may apply \ref{redo} to finally get an $S_{em}>0$-island.
\qed

\begin{corollary}\label{rr}  \, \emph{\textbf{(Reduction of Theorem 1 to 2)}}\label{pks1a} \,   Given a counterexample  $(M^n,g,h)$ to the space-time positive mass conjecture,  we also get an $S_{em}>0$-island.
\end{corollary}

\setcounter{section}{3}
\renewcommand{\thesubsection}{\thesection}
\subsection{MOTS in Compact Initial Data Sets} \label{smq}
\bigskip

We employ the existence of an $S_{em}$-island to build a compact initial data set $\Q^n$ encoding the essential information, carried by a counterexample to Theorem 1, into  a particular, largely torus like geometry of $\Q^n$. This disentangles the technicalities arising from the non-compact asymptotically flat ends from the actual core problems related to the presence of singularities of MOTS. Also it allows the use of skin structural arguments on  $\Q^n$, not available for asymptotically flat spaces.

\subsubsection{Compact Models} \label{ami21}
\bigskip

We start from the unit cube $C^n = [ -1 , 1] \times ...  \times [ -1 , 1] \subset \R^n$ and identify opposite sides of this cube to get the flat torus \[T^n= C^n/\sim=(S^1 \times ... S^1, g_{S^1} \times  ... \times  g_{S^1}),\]
where $g_{S^1}$ is normalized to $length(S^1)=2$. We define the points $1_T:=(1,...,1)/\sim$ and $0_T:=(0,...0)/\sim$ and the
  $(n-1)$-dimensional torus passing through $x \in T^n$:
  \begin{equation}\label{t}
  T^{n-1}[x]:=\{x_1\} \times S^1 \times ... S^1 \subset T^n.
  \end{equation}

Now we use the existence of an \emph{$S_{em}$-island }$(\I^n, g_\I,h_\I)$   with:
 \begin{itemize}
   \item $S_{em}(g_\I,h_\I) >0$ on an  open set, the \emph{core} $U \subset \I^n,\,U \n$, with compact closure.
   \item $(\I^n \setminus U, g_\I,h_\I) \equiv (\R^n \setminus  B_1(0), g_{Eucl},0)$.
 \end{itemize}
We scale $(\I^n, g_\I)$ by some $\gamma \in (0,1/100)$, and get$(\I^n, \gamma^2 \cdot g_\I)$. Now, we insert $(U, \gamma^2 \cdot g_\I) \subset (\I^n, \gamma^2 \cdot g_\I)$ into $T^n$: we delete the ball $B_\gamma(0)$ from $T^n$   and replace it for $U$. This gives a smooth initial data set $(\Q^n_\gamma, g_{n,\gamma},h_{n,\gamma})$ with
\begin{enumerate}
  \item $\Q^n_\gamma$ contains  a flat torus component $\Q T_\gamma^n:=\Q^n_\gamma \setminus U=T^n\setminus  B_\gamma(0)$.
  \item $S_{em}(g_{n,\gamma}) = 0$, since $h \equiv 0$,  on $\Q T_\gamma^n$ and $S_{em}(g_{n,\gamma},h_{n,\gamma})=S_{em}(g_\I,h_\I)> 0$ on  $U$.
\end{enumerate}

Any two parallel tori $T^{n-1}[a_i]  \subset interior(\Q T_\gamma^n)$,  for some $a_1 \neq a_2 \in S^1$,  separate $\Q^n_\gamma$ in two components
\begin{equation}\label{qq}
\Q(flat\,|\,a_1,a_2) \, \dot\cup \, \Q(core\,|\,a_1,a_2):=\Q^n_\gamma \setminus (T^{n-1}[a_1] \cup T^{n-1}[a_2] ),
\end{equation}
where $\Q(flat)$ is just a product of a $T^{n-1}$ with an interval. It is flat relative the metric $g_{n,\gamma}$ and here we also have $h_{n,\gamma} \equiv 0$. The non-trivial part is $\Q(core)$.
It contains the core with $S_{em}(g_{n,\gamma},h_{n,\gamma})> 0$ on  $U$.

\begin{remark} \, From smoothing arguments we can assume that  $\Q^n_\gamma$ is  $C^\infty$-regular. Henceforth, we also assume that $\Q^n_\gamma$ is \emph{orientable}. The non-orientable case can easily be reduced to this one since we always have a orientable double cover. In this cover we can focus on one of the flat torus components and consider the other one as a part of a new and larger core.  \qed\end{remark}

A slight perturbation gives us a geometry on $\Q^n_\gamma$ we choose as the actual starting point for our further discussion:
\begin{proposition} \emph{\textbf{(Compact Model)}}\label{it}\,  For any given $\ve, \gamma >0$, and $l \in \Z^{\ge 4}$, and any two parallel tori $T^{n-1}[a_i]  \subset interior(\Q T_\gamma^n)$,  for some $a_1 \neq a_2 \in S^1$, relative to its flat metric, we can conformally deform $g_{n,\gamma}$ on $(\Q^n_\gamma,g_{n,\gamma},h_{n,\gamma})$ into a smooth metric
\[g_n(\ve,\gamma,l\,|\, a_1,a_2)= v_{\ve,\gamma,l\,|\,a_1,a_2}^{4/n-2} \cdot g_{n,\gamma}, \mm{  with } \]
\begin{enumerate}
  \item $ |g_n(\ve,\gamma,l\,|\, a_1,a_2)-g_{n,\gamma}|_{C^l(\Q^n_\gamma)}\le \ve$,
  \item $S_{em} ( v_{\ve,\gamma,l\,|\,a_1,a_2}^{4/n-2} \cdot g_{n,\gamma},  v_{\ve,\gamma,l\,|\,a_1,a_2}^{2/n-2} \cdot h_{n,\gamma})> 0$ on $\Q^n_\gamma$,
 \item  Both tori $T^{n-1}[a_i] \subset \p \Q(core)$, $i=1,2$, are \textbf{strictly mean convex} relative the new geometry $(\Q(core),g_n(\ve,\gamma,l\,|\, a_1,a_2))$.
\end{enumerate}
\end{proposition}
That is, these tori become barriers for
minimal hypersurfaces, homologous to these tori, to leave the domain $\Q(core)$. Also, we observe that, in (ii) we have, more specifically, $S>0$ and $h \equiv 0$ on $\Q T_\gamma^n$.\\

\textbf{Proof} \quad There are flat tubes $V_i  \subset interior(\Q T_\gamma^n)$ of radius $10 \cdot \zeta >0$, for some $\zeta >0$, around the $T^{n-1}[a_i]$, $i=1,2$, with $V_1 \cap V_2 \v$.
For $x \in V_i $, let $r_i(x)$ denote the \emph{signed distance} of $x$ to $T^{n-1}[a_i]$, so that the sign is positive iff $x \in \Q(core) \cap V_i$.\\

Next, for $f [d,s] (r_i) := s \cdot exp (-d/(r_i+5 \cdot \zeta)$, smoothly extended by zero for smaller $r_i$, and some $\chi \in C^\infty(\R,[0,1])$, with $\chi \equiv 1$ on $\R^{\le 5 \cdot \zeta}$ and
 $\chi \equiv 0$ on $\R^{\ge 6 \cdot \zeta}$ we set
\begin{equation}\label{l}
\phi(d,s)(x):=\chi(r_i(x)) \cdot f [d,s] (r_i(x)).
\end{equation}

Then, from (\ref{fee}), \ref{jesc} and $h \equiv 0$, we can find a large $d \gg 1$ so that for any $s \in (0,1)$:
\[S_{em}((1-\phi(d,s))^{4/n-2} \cdot g_{n,\gamma},0)(x)>0 \mm{ when } |r_i(x)|< 5 \cdot \zeta.\]
Thus on $\Q(flat)$ we have $S_{em} \ge 0$ and we can apply \ref{jesc2} to $U^*=\{x \,|\, |r_i(x)|> 2 \cdot \eta\} \cap \Q(flat)$ and get $S_{em} > 0$ on $\Q(flat)$ without changing $(1-\phi(d,s))^{4/n-2} \cdot g_{n,\gamma}$  on $\{x \,|\, |r_i(x)|\le 2 \cdot \zeta\}$.\\

In turn, to treat $\Q(core)$, we note that we have some (fixed) reservoir of $S_{em}>0$ in the core $U$ and this is essential to compensate the negative $S_{em}$ contribution we unavoidably get from the cut-off by $\chi$ in (\ref{l}. Indeed, we can choose $s>0$ small enough to apply \ref{jesc2}, this time taking $U^*=\{x \,|\, r_i(x)> 2 \cdot \zeta\} \subset \Q(core)$. We get, for small $s>0$,  $S_{em} > 0$ on $\Q(core)$ again without changing $(1-\phi(d,s))^{4/n-2} \cdot g_{n,\gamma}$ on $\{x \,|\, r_i(x)\le 2 \cdot \zeta\}$.\\

At this point, we have  $S_{em} > 0$-geometry on $\Q^n_\gamma$ with some specific shape near the two tori. Since we can choose $s>0$ arbitrarily small, we get both (i) and (ii).\\

  Now, for claim (iii), we recall that the second fundamental form $A_Y(g)$ of a hypersurface  $Y$ with respect to some metric $g$ transforms under conformal deformations to $u^{4/n-2} \cdot g$, for some $C^2$-function $u>0$, according to the formula, cf. [Be],1.163,p.60:
 \begin{equation}\label{mmm}
 A_Y(u^{4/n-2} \cdot g)(v,w) = A_Y(g)(v,w) - \frac{2}{n-2}\cdot {\cal{N}} (\nabla u / u) \cdot g(v,w),
 \end{equation}
where $ {\cal{N}} (\nabla u / u)$ is the normal component of $\nabla u / u$ with respect to $Y$. In our case, we have $A_{T^{n-1}[a_i]}(g_{n,\gamma})=0$, since the tori $ T^{n-1}[a_i]$ are initially totally geodesic, and  $u=1-\phi(d,s)$. This means that $\nabla u$ is a positive multiple of the \emph{outer} normal of $T^{n-1}[a_i] \subset \p \Q(core)$, relative $\Q(core)$.\\

Hence, for $(1-\phi(d,s))^{4/n-2} \cdot g_{n,\gamma}$  on $\{x \,|\, |r_i(x)|\le 2 \cdot \zeta\}$,  the new mean curvature vector, is a negative multiple of the outer normal of $T^{n-1}[a_i] \subset \p \Q(core)$. In other words, the tori $T^{n-1}[a_i]$ are strictly mean convex, relative $\Q(core)$.
\qed

\subsubsection{Almost Minimizers and Stable MOTS} \label{staa}
\bigskip

Here we review the concept of almost minimizers and of stable MOTS.\\

\textbf{Partial Regularity} \,  In the classical treatment  of area minimizing hypersurfaces, as pioneered by De Giorgi, the hypersurface $H$ is considered as the boundary of an open set in $\R^m$.
For these minimizers one only can accomplish a \emph{partial} regularity showing that the hypersurface is largely smooth up to some
singular set  $\Sigma$. The set $\Sigma$ comes without any a priori internal structure, near $\Sigma$ the hypersurface cannot be approximated by hyperplanes, but only by singular cones and $H$ strongly degenerates, in a potentially non-uniform manner, towards $\Sigma$.\\

On the other hand, these methods are remarkably robust. They cover the more general case of \emph{almost minimizers}. This is, despite its dull name, a versatile concept developed, in particular, by Tamanini [T1], [T2],  Bombieri [Bo] and Allard [A]. The following result essentially is [T1],Theorem 1, cf.[L2],Appendix for more details.
\begin{proposition}\emph{\textbf{(Almost Minimizers)}} \,  \label{arm} \, Let $\Omega \subset \R^n$ be open, $E \subset \R^n$ be a Caccioppoli set, which is almost minimizing in $\Omega$ in the sense that
the following estimate, an  almost optimal \textbf{isoperimetric inequality}, holds for some $K > 0$, $\alpha \in (0,1)$
{\small
\begin{equation}\label{iso}
\int_{B_\rho(x)}|D \chi_{E}| - \inf \left\{\int_{B_\rho(x)}|D \chi_{F}|\, \Big| \, F
\Delta E \subset \! \subset B_\rho(x) \right\} \le K \cdot \rho^{n-1+2 \cdot \alpha}
\end{equation}}
for any $x \in \Omega$, $\rho \in (0,R)$, for some $R \in (0,dist(x,\p \Omega))$, where
\begin{itemize}
\item $\chi_A$ is the characteristic function of the set $A \subset \R^n$
\item $\int_\Omega|D \chi_A| := \sup \{\int_\Omega \chi_A \cdot \mbox{div} g  \, d\mu \;| \; g \in C_0^1(\Omega,\R^n), |g|_{C^0} \le 1 \}$ can be thought as the area of $\p A$.
\item $F \Delta E := F \setminus E \cup  E \setminus F$.
\end{itemize}

Then $\p E \cap \Omega$ can be written as a $C^{1, \alpha}$- hypersurface except for some singular set of Hausdorff-dimension $\le n-8$. Also, under blow-ups, that is, under infinite scalings,
we observe a subconvergence to minimal boundaries.
\end{proposition}

\begin{remark} \, 1. In the more general case, where the ambient space is a smooth Riemannian manifold $N^m$, the definition of an \textbf{almost minimizer} in $N^m$ is that of a hypersurface $H^{n-1} \subset N^n$, so that for each point $p \in H$ there is a ball $B \subset N^n$, centered around $p$, so that under some diffeomorphism $\phi$ to $B_1(0) \subset \R^n$, with  $\phi(p)=0$
the condition (\ref{iso}) holds for $\phi(H \cap B)$ near $0$.\\

2. Area minimizing hypersurface clearly are almost minimizer. In this case, the left hand side of (\ref{iso}) vanishes. But  (\ref{iso}) also means that any \emph{smooth} hypersurface in a compact manifold is almost minimizing. In other words, the condition really characterizes the highly curved and singular portion of the hypersurface as asymptotically simulating proper area minimizers.  \qed
\end{remark}

\textbf{MOTS} \, The concepts and results concerning MOTS in this section have been developed in a series of papers by Andersson, Eichmair, Galloway, Mars, Metzger, Schoen and Simon [AMS1-2], [AM], [E2], [EM] and [GS]. An instructive survey is given in [M]. For the purposes of this paper we use almost minimizers already in the definition of MOTS, to have a dense open set where the various operators and curvatures are properly defined.

\begin{definition}\emph{\textbf{(Stable MOTS)}} \, \label{mor} Let $(M^n,g,h)$ be an orientable and compact initial data set and $H^{n-1}\subset M^n$ a compact almost minimizer.\\

 $H^{n-1}$ is called a \textbf{marginally outer trapped surface}, abbreviated \textbf{MOTS}, if
\begin{equation}\label{mots}
\theta^+_H:=tr_H h_{M|H} +  tr_H h_H=0
\end{equation}
 where $tr_M h$ denotes the trace of $h$ only over the tangent space of $H$.\\

$H^{n-1}$  is called a \textbf{stable MOTS} if it is MOTS so that there exists a smooth function $f>0$ on $H^{n-1} \setminus \Sigma_H$, such that $L_H^{MOTS} f \geq 0$, where, for any smooth $v$:
\begin{equation}\label{lh}
L_H^{MOTS} v :=-\Delta_H v +2 \langle X, \nabla v\rangle  + (div_H X - |X|^2+Q_H) \cdot v,
\end{equation}
where $Q_H$ is the container for the curvature terms
\begin{equation}\label{lh2}
 \textstyle Q_H=\frac{1}{2}S_H-\mu-J(\nu)-\frac{1}{2}|h_{M|H}+h_H|^2,
\end{equation}
\end{definition}
The terms which appear in (\ref{lh}) and  (\ref{lh2}) are defined as follows:
 \begin{itemize}
 \item $\nu$ is the outer unit normal field on $H \setminus \Sigma_H$, locally viewed as a boundary of some open set. The choice of what is inside or outside is made to match the choice of the normal direction in (\ref{mots}),
   \item $X$ is the tangential part of the vector field dual to $h(\nu, \cdot)$ on $H$,
   \item $h_{M|H}$ is the restriction of $h$ on the ambient space $M$ to vectors tangent to $H$,
   \item $h_H$ is the second fundamental form of $H \subset M$,\\
\end{itemize}
A comparison of $L_H^{MOTS}$ with the operator $L_H^Q=-\Delta_H +Q_H$ gives the following notable result
\begin{lemma}\label{smo}
Let $H$ be a stable MOTS in an initial data set $(M,g,h)$. Then, for every $C^1$ function $v$ compactly supported on $H \setminus \Sigma$, we have
\begin{equation}
\int_H |\nabla v|^2  + Q_H \cdot  v^2  \,dA \geq 0.
\end{equation}
\end{lemma}
\textbf{Proof} \quad  This result was proved in [GS], cf. also [M],Ch.3.1 for a nice presentation.\qed

Now we show that $\Q^n$ carries such a stable MOTS homologous to  the tori $T^{n-1}[x]$ and observe that they contain a nearly flat and regular torus component. Here we think of $S^1=[-1,1]/\sim$ and simply write $\eta \in \R$, for $|\eta|<1$, for $\{\eta\}/\sim \, \in S^1$:

\begin{proposition}\emph{\textbf{(Existence of Stable MOTS in $\Q^n$)}}   \label{exmo}   For any $\ve^* >0$, $\gamma^* \in (0,1/100)$ and $l^* \in \Z^{\ge 4}$, there are some $\ve \in (0,\ve^*)$, $\gamma \in (0,\gamma^*)$ and   $l = l^*+1$, so that $(\Q^n_\gamma,g_n(\ve,\gamma,l))$, as in \ref{it},  contains a \textbf{stable MOTS} $H(\ve,\gamma,l | -2  \cdot \gamma, 2  \cdot \gamma) \subset \Q^n_\gamma$ with
\begin{itemize}
  \item $H \subset \Q(core \, | -2  \cdot \gamma,2  \cdot \gamma)$ and $H$ is homologous to $T^{n-1}[a_1]$,
  \item $H \cap  \Q T^n_{\gamma^*}$ is a $C^{l^*}$-regular hypersurface.
  \item  $H \cap  \Q T^n_{\gamma^*}$  is almost isometric to $T^{n-1}[0_T] \setminus B_{\gamma^*}(0_T)$:
  \begin{equation}\label{f0}
 H \cap   \Q T^n_{\gamma^*}= exp_\nu(\Gamma)  \cap \Q T^n_{\gamma^*},
 \end{equation}
for  the image $exp_\nu(\Gamma)$ of some $C^{l^*}$-regular section $\Gamma$ of the normal bundle $\nu$ over $T^{n-1}[0_T] \setminus   B_{\gamma^*}(0_T)$,  under the exponential map $exp_\nu$ of $\nu$, with
 \begin{equation}\label{f}
 |\Gamma|_{C^{l^*}\left(T^{n-1}[0_T] \setminus   B_{\gamma^*}(0_T)\right)} <\ve^*.
 \end{equation}
\end{itemize}
\end{proposition}

\textbf{Proof} \quad  The existence result is a direct combination of [AM],Th.1.1, alternatively [E2],Th.1.1, and [EM],Th.2.3. Some partial results in these reference are stated only for the smooth case of stable MOTS in dimension $\le 7$, but what we need here can be readily extended to the general cases, specifying the result as a potentially singular almost minimizer.\\

The key geometric property we exploit to apply these results is \ref{it}(iii). It ensures that the MOTS, which is a stable minimal hypersurface outside the core,  remains in $\Q(core \, | -2  \cdot \gamma,2  \cdot \gamma)$. This is the counterpart to [L1], Prop.2.1 in in the case of area minimizing hypersurfaces, albeit, the techniques to show the existence of these MOTS follow different lines. The other statements are transcripts from the corresponding  Prop.2.3 in [L1]. They remain valid since Allard regularity theory equally applies to almost minimizers as in the argument for [L1],Prop.2.3. \qed

\subsubsection{Skin Structures and Surgeries} \label{uom}
\bigskip

In this final section we gather  the conclusions from the previous chapters, but also we invest major parts of [L1]-[L4], to show that the potentially singular stable MOTS we have obtained in \ref{exmo} can be deformed to torus-like spaces with positively mean curved boundaries. But an inductive argument, explained in [L1], shows that such spaces cannot really exist and this contradiction
finishes the argument for our Theorems.\\

\textbf{Skin Structures in a Nutshell} \, We consider a connected almost minimizer $H$ with singular set $\Sigma$. Typically we assume that $\Sigma \n$. The theory also applies when $H$ smooth, but then many statements  become  trivial.\\

Under blow-ups, that is, infinite scalings around singular points, almost minimizers flat norm subconverge to minimal boundaries in $\R^n$. This shows that the skin structural techniques developed in [L2]-[L4] equally apply to almost minimizers since all estimates are actually inherited from the limiting spaces, which, in turn, form a compact ensemble of model spaces.\\

The basic notion is that of a \textbf{skin transform}  $\bp_H$. It  is a non-negative measurable function naturally defined on $H \setminus \Sigma$: it
commutes with the compact convergence of sequences of compact almost minimizers to area minimizers, and of  Euclidean oriented boundaries to other such boundaries, and it satisfies the following axioms
\begin{itemize}
    \item $\bp_H \ge |h_H|$ and for any $f \in C^\infty(H \setminus \Sigma,\R)$ compactly supported in $H \setminus \Sigma$ and some $\tau = \tau(\bp,H) \in (0,1)$ we have the Hardy type inequality
\begin{equation}\label{hi}
\int_H|\nabla f|^2  + |h_H|^2 \cdot f^2 dA \ge \tau \cdot \int_H \bp_H^2\cdot f^2 dA.
\end{equation}
    \item $\bp_H \equiv 0$, if $H \subset M$ is totally geodesic, otherwise, $\bp_H$ is  strictly positive.
    \item When $H$ is not totally geodesic, we define $\delta_{\bp}:=1/\bp$, the \emph{${\bp}$-distance}. It is $L_{\bp}$-Lipschitz regular, for some constant
        $L_{\bp}=L(\bp,n)>0$:
        \[|\delta_{\bp}(p)- \delta_{\bp}(q)|   \le L_{\bp} \cdot d(p,q), \mm{ for } p,q \in  H \setminus \Sigma. \]
\end{itemize}
In [L2] we described procedures to define such a skin transform which seamlessly extend to almost minimizers. Here, and henceforth, $h_H$ is the form we get from the embedding of $H$ in some existing  ambient space. That is, as in \ref{me}, $h_H$ transforms, under conformal deformations, into the new form without the correction term of (\ref{h01}).\\

For any compact and connected almost minimizer $H$ we find a quantitative sort of connectedness for $H \setminus \Sigma$. It is a \textbf{skin uniform} space: for any two $p,q \in H$, there is a rectifiable path $\gamma: [a,b] \ra H$, for some $a <b$, with $\gamma(( a,b))\subset H \setminus \Sigma$,
so that for any given skin transform $\bp$, there is some $s_H \ge 1$ with
\[l(\gamma)  \le s_H \cdot d_{g_H}(p,q) \mm{ and }  l_{min}(\gamma(z)) \le s_H \cdot \delta_{\bp}(z),\mm{ for any } z \in \gamma_{p,q}.\]

To check this for almost minimizers, we adapt the argument from [L2]. For the local connectedness of $H \setminus \Sigma_H$, we observe that after scalings $H$ can locally be approximated by an oriented minimal boundary $H^* \subset \R^n$. We  transfer the estimates to derive the local connectedness of $H^*$ in [L2], Lemma 4.12 to $H$. The skin uniformity then follows from the same blow-up and network argument as in [L2],Ch.4.3.\\

The skin uniformity is the main prerequisite to get \textbf{hyperbolic unfoldings} of $H \setminus \Sigma$ to a conformally equivalent, complete Gromov hyperbolic space with bounded geometry, cf. [L3],Th.2. As was shown by Ancona [An], the potential theory of uniformly elliptic operators (satisfying a certain coercivity condition) on such spaces is easy to describe. The point is that this understanding can be transferred back to the original space $H \setminus \Sigma$. (In this context we no longer use minimality properties of $H \setminus \Sigma$.)\\

This strategy can be applied to \textbf{skin adapted operators}, that is, linear second order elliptic operators $L$, so that
\emph{\begin{enumerate}
     \item $L$ does not degenerate faster than $\bp^2$ when we approach $\Sigma$.
     \item $L$ satisfies a weak coercivity condition: we require that there is a supersolution $s
         >0$ of $L \, u=0$ and some $\ve >0$, so that: $L \, s \ge \ve \cdot \bp^2 \cdot s.$
\end{enumerate}}
What we get are controls for the asymptotic behavior of solutions $u>0$ of $L \, w=0$, when we approach $\Sigma$, similar to those we otherwise merely expected for the Laplacian on a smoothly bounded
Euclidean domain. \\

\textbf{Conformal Laplacians} \,  The analytic key result we establish now is that on stable MOTS in an initial data set, satisfying the DEC, the conformal Laplacian is skin adapted. To this end we start with a Hardy inequality for stable MOTS extending the Hardy inequality  (\ref{hi}) to the case of initial data sets.

\begin{proposition}\emph{\textbf{(Hardy Inequality for Stable MOTS)}}\label{lh0} \, Let $H$ be a stable MOTS in a compact initial data set $(M,g,h)$. Then  there is a constant $a_H>0$, so that for  any smooth function $v$ with $supp \: v \subset H \setminus \Sigma$, we get
{\small \begin{equation} \label{wo} \int_H |\nabla v|^2  + |h_{M|H}+h_H|^2\cdot  v^2  \,  dA \ge a_H \cdot \int_H \bp^2\cdot  v^2  \,  dA
\end{equation}}
\end{proposition}

\textbf{Proof} \, We use a particular model of a skin transform to derive (\ref{wo}) as a initial data set variant of (\ref{hi}). The skin transform $\bp _1$ of [L2],Th.1 and Def.2.7 is given by
\begin{equation} \label{wen} \bp _1(x):= \sup \{c \,| \, x \in \overline{\U_c}\}.\end{equation}
for any $x \in H \setminus \Sigma$, where $\U_c$ is the $1/c$-distance tube of $|h_H|^{-1}([c,\infty))$. For totally geodesic $H$, we set $\bp_1\equiv 0$. Then we
argue as in the proof of [L2],Prop.3.3 and 3.4  where we established the standard version (\ref{hi}).\\

To this end, we start with a narrow neighborhood $W$ of $\Sigma$ and neighborhood $V$ of $H \setminus W$, compactly supported  in $H \setminus \Sigma$. We cover $W$ by a locally finite collection of balls $B_{\mu /\bp_1(p_i)}(p_i)$, for suitable $p_i$,$i \in \Z^{\ge 1}$ and $\mu \in (0,1/2)$, with uniformly upper bounded intersection numbers, a skin adapted cover, cf.[L2], Ch.3.1. The strategy of [L2],Prop.3.3 and 3.4,  we recycle here, is to get (\ref{wo}) from a common lower bound of the $\bp_1$-weighted Neumann eigenvalues of  $-\Delta + |h_{M|H}+h_H|^2$ on $V$ and  on the balls $B_{\mu /\bp_1(p_i)}(p_i)$. For suitably large $V$ this is obvious since $|h_{M|H}+h_H|$ diverges towards $\Sigma$.\\

For the balls $B_{\mu /\bp_1(p_i)}(p_i)$ we assume a subconvergence, for $i \ra \infty$:
{\small \begin{equation}\label{sub}
\inf_{f \in C^\infty (B_{\mu /\bp_1(p_i)}(p_i)), f \not\equiv 0}\frac{ \int_{B_{\mu /\bp_1(p_i)}(p_i)} |\nabla f|^2  + |h_{M|H}+h_H|^2\cdot f^2 dA}{\int_{B_{\mu /\bp_1(p_i)}(p_i)} \bp_1^2\cdot f^2 dA} \ra 0
\end{equation}}

As in [L2],Eq.(29) these $\bp_1$-weighted Neumann eigenvalues are scaling invariant.  (\ref{sub}) says that also the unweighted first Neumann eigenvalues of  $-\Delta + |h_{M|H}+h_H|^2$ and a suitably $L^2$-normalized first eigenfunction on $\bp_1(p_i) \cdot B_{\mu /\bp_1(p_i)}(p_i)$ subconverge to zero and some associated positive eigenfunction. But these scaled ball smoothly subconverge to a ball $B$ with $h_H \not\equiv 0$ in an oriented Euclidean minimal boundary.\\

$|h_{M|H}|$ remains bounded on $H$,  while  $ |h_H|$ diverges towards $\Sigma_H$. Thus we get, on $\bp_1(p_i) \cdot B_{\mu /\bp_1(p_i)}(p_i)$, a $C^3$-convergence of $|h_{M|H}|$ to zero. While $\inf  h_H$ remains positively lower bounded. Then we get, as in [L2],Lemma 3.8, from the definition (\ref{wen}) of $\bp_1$,  that the first Neumann eigenvalues of  $-\Delta + |h_H|^2$ on this ball $B$ must vanish, contradicting $h_H \not\equiv 0$. \qed

The conformal Laplacian $L_H:=-\Delta  +\frac{n-2}{4 (n-1)} \cdot scal_H$ appears in the scalar curvature  transformation rule under conformal deformations $g_H \mapsto u^{4/n-2} \cdot g_H$, we recall from Ch.\ref{ami0}.A (\ref{cvd}):
\begin{equation}\label{stl}
\textstyle S(u^{4/n-2} \cdot g_H) \cdot u^{\frac{n+2}{n-2}} = \gamma_n \cdot L_H\,  u, \mm{ for } \gamma_n = \frac{4 (n-1)}{n-2},
\end{equation}
valid for any $C^2$-regular function $u>0$. This makes an understanding of the analysis of $L_H$ valuable in scalar curvature geometry.

\begin{proposition}\emph{\textbf{($L_H$ on Stable MOTS)}}\label{lh0} \, Let $H$ be a stable MOTS in a compact initial data set $(M,g,h)$ with $S_{em} \ge 0$. Then, $L_H$ is a skin adapted operator.
\end{proposition}
In this case (of a symmetric operator) the weak coercivity requirement (ii) of the skin adaptedness conditions is equivalent to
{\small \begin{equation}\label{cf}
\int_H  f  \cdot  L_H f  \,  dA  \ge \tau \cdot \int_H \bp^2 \cdot f^2   \, d A, \mm{ for some }\tau >0,
\end{equation}}
for any smooth function $f$ on $H$ with $supp \: f \subset H \setminus \Sigma$, cf.[L3],Lemma 5.3.\\

\textbf{Proof} \quad The DEC, that is,  $S_{em} \ge 0$ saying that  $\mu \ge |J|$ implies
{\small \begin{equation}\label{dee}
\frac{1}{2}S_H-\frac{1}{2}|h_{M|H}+h_H|^2 \ge Q_H=\frac{1}{2}S_H-\mu-J(\nu)-\frac{1}{2}|h_{M|H}+h_H|^2
\end{equation}}
Thus, for any smooth function $v$ with $supp \: v \subset H \setminus \Sigma$, we get from \ref{smo}
{\small \begin{equation}\label{r}
\int_H |\nabla v|^2  + \frac{1}{2} (S_H-|h_{M|H}+h_H|^2) \cdot  v^2  \,  dA  \ge \int_H |\nabla v|^2  + Q_H \cdot  v^2  \,dA \ge 0
\end{equation}}
Also directly from the definition of $L_H$:
{\small \begin{equation}\label{sim}
\int_H  v  \cdot  L_H v  \,  dA  \ge \frac{n-2}{4 (n-1)} \cdot \int_H | \nabla v|^2  d A + \int_H \frac{n-2}{2 (n-1)} \cdot | \nabla v|^2 + \frac{n-2}{4 (n-1)} \cdot S_H  \cdot  v^2 \, d A
\end{equation}}
Combing the two inequalities (\ref{r}) and (\ref{sim}) we have
{\small  \begin{equation}\label{wer}
 \int_H  v  \cdot  L_H v  \,  dA  \ge \frac{n-2}{4 (n-1)} \cdot \int_H |\nabla v|^2  + |h_{M|H}+h_H|^2\cdot  v^2  \,  dA
\end{equation}}
Finally, we combine (\ref{wer}) and (\ref{wo}) to infer the validity of (\ref{cf}). \qed

At this point we can jump into the routine in the proof of the Riemannian positive mass conjecture, in [L1], where in place of the MOTS of \ref{exmo} in an ambience with $S_{em}>0$, we had a similar area minimizer in a  $S>0$-ambience. In both cases the conformal Laplacian is skin adpated. We can apply the identical surgery techniques, explained in [L1],Ch.2.3 and 2.4 and get the following spaces $N^{n-1}$, from  deformations and surgeries applied to the  MOTS $H^{n-1} \subset (\Q^n_\gamma,g_n(\ve,\gamma,l))$ of \ref{exmo}:

\begin{proposition}\emph{\textbf{(MOTS Surgery)}}\label{id}  We assume the spaces $\Q^n_\gamma$ described in \ref{it} exist. Then, for any $\ve>0$, $\gamma \in (0,1/100)$, $l \in \Z^{\ge 4}$, there is a $C^l$-regular $n-1$-dimensional Riemannian manifold $(N^{n-1}(\ve,\gamma,l),g_{n-1}(\ve,\gamma,l))$ with the following properties
\begin{enumerate}
\item $(N^{n-1},g_m)$ is a compact $S>0$-manifold with boundary $\p N^{n-1}$.
\item The boundary $\p N^{n-1}$ may be empty and for $n-1 <7$ it is empty.
\item  $ N^{n-1}$ contains a torus component $NT^{n-1}$, nearly isometric to $T^{n-1}\setminus B_{\gamma}(0)$:
\[\mm{There is a $C^l$-regular diffeomorphism } F_{n-1}(\ve,\gamma,l): NT^{n-1} \ra T^{n-1}\setminus  B_{\gamma}(0)\mm{ with }\] \[ |F_{n-1}^*(g_{flat})- g_{n-1}|_{C^l(NT^{n-1})} \le \ve.\]
\item  When $\p N^{n-1} \n$, then $\p N^{n-1}  \subset N^{n-1} \setminus NT^{n-1}$ is  $C^3$-smooth and compact with positive mean curvature.
\end{enumerate}
\end{proposition}

Now we can appeal to [L1],Prop 1.2 and Cor.1.4. They show that there cannot be such families of manifolds $N^{n-1}$. This contradiction proves that the $\Q^n_\gamma$ of \ref{it},
the $S_{em}$-islands and, finally, any sort of counterexamples to the space-time positive mass conjecture, cannot exist. \qed

\footnotesize
\renewcommand{\refname}{\fontsize{14}{0}\selectfont

\end{document}